\documentclass[11pt]{amsart}
\usepackage{amsmath, amsthm, amsfonts,amssymb} 

\theoremstyle{plain}
\newtheorem{theorem}{Theorem}[section]
\newtheorem{lemma}[theorem]{Lemma}
\newtheorem{proposition}[theorem]{Proposition}
\newtheorem{corollary}[theorem]{Corollary}

\theoremstyle{definition}
\newtheorem{definition}[theorem]{Definition} 
\newtheorem{validrule}{Rule}

\theoremstyle{remark}

\newtheorem*{remark}{Remark}

\newcommand{\imprho}{\stackrel{\rho}{ \Rightarrow }}
\newcommand{\notrho}{\stackrel{\rho}{\neg}\!}
\newcommand{\truesent}{{\mathcal{T}}}
\newcommand{\falsesent}{{\mathcal{F}}}
\newcommand{\bN}{{\mathbb{N}}}
\newcommand{\turnrho}{{\,\vdash_{\!\!\!\rho}\,}}
\newcommand{\dimp}{{\rightarrow }}
\newcommand{\iffrho}{\stackrel{\rho}{ \iff }}
\newcommand{\ahalt}{{\mathcal{H}}}
\newcommand{\cB}{{\mathcal{B}}}
\newcommand{\imprhoone}{\stackrel{\rho_1}{ \Rightarrow}}
\newcommand{\imprhotwo}{\stackrel{\rho_2}{ \Rightarrow }}

\newcommand{\goldbach}{\textsc{Goldbach}}
\newcommand{\halt}{\textsc{halt}}
\newcommand{\truealg}{\textsc{true}}
\newcommand{\proverho}{\textsc{prove}_\rho}
\newcommand{\identity}{\textsc{identity}}
\newcommand{\deduce}{\textsc{deduce}}
\newcommand{\asand}{\textsc{and}}
\newcommand{\asor}{\textsc{or}}
\newcommand{\disjunctelim}{\textsc{d-elim}}
\newcommand{\sneg}{\textsc{s-neg}}
\newcommand{\mpone}{\textsc{mp$(\rho_1)$}}
\newcommand{\denyc}{\textsc{deny$(C)$}}
\newcommand{\curry}{\textsc{Curry}}

\begin{document} 

\title{Stability and Paradox in Algorithmic Logic}

\author{Wayne Aitken, Jeffrey A.~Barrett}

\date{November 20, 2005}

\address{Cal.~State, San Marcos, CA 92096, USA}
\email{waitken@csusm.edu}

\address{UC Irvine, Irvine, CA 92697, USA}
\email{jabarret@uci.edu}


\begin{abstract}
There is significant current interest in 
type-free systems
that allow flexible self-application.
Such systems are of interest in property theory,
natural language semantics, 
the theory of truth, 
theoretical computer science, 
the theory of classes, 
and category theory.
While there are a large variety of proposed type-free systems,
there is a particularly natural 
type-free system that we believe is prototypical:
the logic of recursive algorithms.
\emph{Algorithmic logic} is the
study of basic statements concerning 
algorithms and the algorithmic rules
of inference between such statements. 
As shown in~\cite{firstpaper04},
the threat of paradoxes, such as the Curry paradox,
requires care in implementing
rules of inference in this context. 
As in any type-free logic, some traditional rules will fail.
The first part of the paper develops
a rich collection of rules of inference which do
not lead to paradox. The second part identifies traditional rules
of logic which are paradoxical in algorithmic logic, and so 
should be viewed with suspicion in type-free logic.
\end{abstract}
\maketitle


\section{Introduction}\label{s1}

In second-order logic, one distinguishes between two types of objects.
\emph{First-order objects} are the basic objects of interest. \emph{Second-order
objects} are the properties and classes, the functions and operators for the
first-order objects. There are, however, situations in which this division
is unnatural. When one wants, for whatever purpose,
to mix the first-order and second-order universes, one is reminded of the
reasons for their original separation: the paradoxes.

This paper is a study in \emph{type-free logic}.
The goal of type-free logic is to find consistent, natural, and flexible ways to handle
type-free systems where the second-order objects are not
separated from the first-order objects, but are in some sense part of
the first-order universe. 
In type-free logic one desires enough flexibility to including meaningful self-application 
and self-containment: functions and properties should be able to apply to themselves
and collections should be able to contain themselves. 
One also wants natural or ``na\"ive'' comprehension and functional abstraction
principles to create collections and functions. One might also want a truth predicate.

Standard ZFC set theory fails these desiderata.
Functions may be applied to functions, and collections may contain collections,
but within limits:
a function cannot be a member of its domain and a collection cannot contain itself.
In ZFC only a well-behaved, well-founded part of the second-order universe is
allowed inside the first-order universe. In fact, due to well-foundedness,
the set-theoretic universe can be regarded as a typed theory
where the universe is typed by ordinal rank.
Meaningful self-application is blocked:
the relation $x\in y$ is automatically false unless $x$ has a lower 
ordinal rank than $y$.
ZFC does not allow a universal set (and even extensions of ZFC
such as GB or MK class theory
do not allow a universal class that contains the universal class).
A typical example of the limitations of ZFC in this regard
concerns the status of the self-composition function $T$:
given a function $f$ with domain and codomain the same class, 
defined $Tf$ to be $f\circ f$. 
This natural function is not an object
in the ZFC universe.\footnote{Of course, \emph{restricted} versions of 
$T$ can be defined 
to apply to functions on a fixed domain.}
Now there is nothing preventing one
from studying $T$ in set theory,
but the point is that it is external to the set-theoretical universe~$V$, a universe
intended to be rich enough for all of mathematics. This example is not
atypical. Set theorists work outside of $V$ whenever
properties about the intersection $\cap$ or union $\cup$ operators are discussed:
the set-theoretical operators $\cap$ or $\cup$ are not themselves objects
of the set theoretical universe~$V$.

Whenever a type-free system is considered with more expressive
power than ZFC set theory, 
for example a system with unrestricted comprehension,
the threat of paradox emerges anew. So some part or another of 
the traditional 
logics must be restricted.
Nevertheless, there is significant current interest in 
type-free systems due to applications in 
property theory,
 natural language semantics, 
the theory of truth, 
theoretical computer science, 
the theory of classes, 
and category theory.
In property theory, it is desirable, indeed arguably essential, 
that every open formula in
a language should determine an associated
object called a \emph{property}.
There should be a truth predicate that behaves in a manner similar to the
truth predicate in natural language.
In class theory there should be nothing preventing a class from
containing itself. In fact, any restriction on class comprehension
seems artificial.\footnote{This is in contrast to \emph{sets} where
good conceptual reasons have been given to restrict the comprehension principle.}
There should be a universal class, and this class should contain itself.
And given its role as an organizing principle of contemporary mathematics,
there should to be a more satisfying way to develop category theory than
by employing the current large/small category distinction.

There are a large variety of proposed type-free systems\footnote{Recent
examples include \cite{field04}, \cite{orilia00}, and \cite{cantini96}.
Recent examples from the substructural tradition
include 
\cite{terui04},
\cite{cantini03},
\cite{girard98}, and \cite{weir98}.
See the bibliographies of these works for earlier examples.
The articles in~\cite{russell04}, especially those by 
A.~Cantini,
S.~Feferman,
H.~Field, H.~Friedman, H.~Sturm, and 
K.~Wehmeier,
show the contemporary interest in 
type-free systems and in
 strong forms of comprehension and abstraction.}
which are provably free from contradictions engendered by
paradoxes, and which restrict the traditional logics in one way or another. 
What remains is the question of which
type-free systems are the most compelling.
One obvious criterion is that a type-free system should not introduce
artificialities worse that the artificiality,
discussed above, of separating first-order objects from second-order objects.

In this paper we introduce a promising methodology
for developing a natural type-free system.
Perhaps the most common strategy 
is to start with some form of a classical logic with
a  na\"ive comprehension principle, then to weaken
it until it is consistent. 
But herein lies the problem: it is unclear what to weaken.
Our strategy, on the other hand,
is to begin with a naturally occurring  type-free system then to
investigate the logical properties it in fact possesses. 
The hope is that this naturally occurring
type-free system will serve as a fruitful model for strong type-free systems
more generally.

Perhaps the most natural type-free system is ordinary language, but
for our purpose it is hopelessly intractable. 
The universe of recursive algorithms, however, is both natural and tractable.
If we fix a framework for algorithmic description, then self-application of
algorithms is possible, indeed commonplace. If we focus on logical
operators that can be defined algorithmically, a rich type-free logical structure
emerges. \emph{Algorithmic logic} is the study of this 
type-free system.

This paper is the second of a series designed to introduce and examine algorithmic logic.
The first~\cite{firstpaper04}, 
a short and informal introduction to the subject,  focussed on 
the challenge of the Curry paradox. 
The Curry paradox is the first test of any type-free system
containing implication. 
The present paper builds on the lessons learned from the
first, but is independent of it. It begins the formalization
and careful study of algorithmic logic.
The main task here 
is to understand which traditional rules of propositional logic are safe
and which are problematic in algorithmic logic. 
A third paper~\cite{thirdpaper06} will discuss
the principle of unrestricted functional abstraction in algorithmic logic.
Fredrick Fitch \cite{fitch69} sought a type-free logic with
an unrestricted abstraction principle; he regarded any
restriction on abstraction
as artificial and undesirable. 
Since algorithmic logic is both type-free and allows for a strong
abstraction principle,
our work can be viewed as part of the Fitch--Curry--Myhill 
tradition.\footnote{We recommend \cite{myhill84}, \cite{feferman84}, 
and \cite{cantini96}
as interesting introductions to type-free logic. 
We have found \cite{restall94} to be a helpful introduction to
the substructural tradition.}


\section{Algorithmic Logic}\label{s2}

The basic objects of algorithmic logic are \emph{algorithmic statements}.
An algorithmic statement is an assertion of the form \emph{algorithm
$\alpha$ with input $u$ halts with output $v$}. 
The assertion $4! = 24$ can be understood as a true algorithmic statement, 
where the algorithm is one designed to calculate the factorial function,
the input is $4$, 
and the output is $24$.

Algorithmic statements can be more subtle.\footnote{In this paper
\emph{algorithms} will be limited to recursive algorithms. 
With this restriction, the collection of
algorithmic statements
can be seen to be in some sense 
equivalent to the collection of $\Sigma_1$-statements
in first-order arithmetic.} 
Consider Goldbach's 
conjecture that every even number greater than
two is the sum of two prime numbers.
The \emph{negation} of Goldbach's conjecture can
be understood as the algorithmic statement that $\goldbach$ 
halts with output $0$ when run with input $0$,
where $\goldbach$ is the algorithm that checks each even number in turn,
beginning with four, and outputs $0$ if it ever finds 
a number that cannot be represented as the sum of two primes.  
Note that if Goldbach's conjecture is true, then the algorithm $\goldbach$ will simply fail to halt regardless of input.  

An algorithmic statement can be false in two ways.  
It can be false because the algorithm halts with an output different from the one specified.  Such statements are \emph{directly false}.  
Or it can be false because the algorithm fails to halt.  
Such statements are \emph{indirectly false}.

The assertion that a specified algorithm \emph{halts} on a specified input 
can also be understood as an algorithmic statement. 
Consider the algorithm $\halt$ that takes as input a pair~$[\alpha, u]$
and runs as a  subprocess the algorithm $\alpha$ with input $u$. 
The algorithm $\halt$
outputs~$1$ if the subprocess halts; otherwise $\halt$ itself does not halt.
So the algorithmic statement asserting that $\halt$ outputs
$1$ on input $[\alpha, u]$ is true if and only if $\alpha$ halts on input $u$.

The algorithm $\halt$ is an example of an \emph{algorithmic predicate},
a predicate that can be represented by an 
algorithm that outputs $1$ if and only if the predicate is true of the input.
 We require that if a predicate algorithm halts at all, it outputs~$0$ or~$1$.
Algorithmic predicates are the basic
\emph{internal} predicates of algorithmic logic.

There is an internal truth predicate $\truealg$ for algorithmic statements.
The algorithm $\truealg$ expects as input a triple $[\alpha, u, v]$ representing an algorithmic statement with
specified algorithm~$\alpha$, specified input~$u$, and specified output~$v$.
First $\truealg$ runs the subprocess $\alpha$ 
with input~$u$. If this subprocess halts with output $v$, then $\truealg$ outputs~$1$.  
If the subprocess halts with output not equal to $v$, then $\truealg$ outputs~$0$.  
If the subprocess fails to halt, then $\truealg$ also fails to halt.

Closely related to the truth predicate, is an algorithmic predicate 
corresponding to directly false.  
Because of the halting problem, however, there is no algorithmic predicate 
corresponding to false:
the \emph{external} property of being false is one that cannot be 
expressed internally.

Finally, \emph{algorithmic connectives} can be defined in terms of algorithmic 
predicates.  The algorithmic conjunction $\wedge$ and disjunction
$\vee$ behave as expected,  
but the algorithmic conditional $\imprho$
requires special care.  
Here the conditional is indexed by a library~$\rho$ of
inference rules.  The algorithmic statement $A \imprho B$ means
that the algorithmic statement $B$ can be deduced from the algorithmic
statement $A$ using the rules in the library $\rho$. 
Because of its definition, the connective
$\imprho$ can be used to define an internal predicate~$\proverho$.
The connective $\imprho$ is also used to define
negation~$\notrho$\, .

Since algorithms can take algorithms as input, 
as in the case of $\halt$ above, algorithmic logic
is inherently self-referential and so is 
essentially type-free. Consequently,
special care must be taken to avoid contradiction:
the rules of the library $\rho$ must be carefully evaluated for
validity.   
Indeed, in an earlier paper~\cite{firstpaper04} we show
 that the rule \emph{modus ponens}
for $\imprho$
cannot be included in a sufficiently rich library $\rho$ 
without rendering
the rule itself invalid.
If \emph{modus ponens} is included in
such a library, an algorithmic version of the Curry
paradox results in a contradiction.

The first part of the present paper introduces rules for algorithmic logic
that form a \emph{stable base}: a collection of rules 
that can be safely extended to form stronger
valid collections. The second part of the paper
presents a list of \emph{paradoxical rules}:  
traditional rules of logic that can be shown to be
invalid when in a sufficiently rich library, usually through
arguments akin to those found in the Russell and Curry paradoxes.


\section{Conventions for Algorithms}\label{s3}

Rather than stipulate a particular theoretical framework for the discussion of 
algorithms, 
we require that a suitable framework behaves as follows.

Anything that can be input into an algorithm is called a \emph{datum}.  
Data include natural numbers and algorithms.  In addition,
if $a_1 \ldots a_k$ are data, the list~$[a_1, \ldots , a_k]$ is itself a datum.
Every algorithm accepts exactly one input datum and either does not halt or 
halts with exactly one output datum.  
If an algorithm requires or produces multiple data, 
the data are packaged in a single input or output list respectively.
Any datum is an allowable input whether or not 
it is consistent with the intended function of the algorithm.
Typically, we will not specify what an algorithm does with 
an unexpected input datum.

Every datum has a positive integer~\emph{size},
and there are only a finite number of data of a given size.
The size of a list is strictly greater than the sum of the sizes of the items of the list.  
A \emph{process} is a pair consisting of an algorithm and an input.   
Every halting process has a positive integer \emph{runtime}.  
If a parent process runs one or more subprocesses in its execution, 
then the runtime of the parent process is strictly greater than the sum of the runtimes of the halting subprocesses.  

An \emph{algorithmic statement} 
can be represented as a datum:
 If $\alpha$ is the specifed algorithm, $u$ the specified input, and $v$ 
the specified output, then the list $[\alpha, u, v]$ represents the corresponding
algorithmic statement.

The identity algorithm $\identity$ simply outputs a copy of its input.
The algorithmic statement $[\identity, 0, 0]$ is denoted~$\truesent$.
This algorithmic statement is true.
Similarly,   the algorithmic statement $[\identity, 0, 1]$ is denoted~$\falsesent$.
This statement is directly false since~$0\neq1$.  


\section{Deduction}\label{s4}

There is an algorithmic predicate for deduction. 
This deduction predicate depends on a library of rules
instantiated by an algorithmic sequence.
\begin{definition}\label{def41}
An \emph{algorithmic sequence} is an algorithm which halts 
for every positive integer input.
If $\alpha$ is an algorithmic sequence, then $\alpha_n$ denotes the output of
$\alpha$ applied to the integer~$n$.  
\end{definition}

Informally, a rule is an algorithm which expects
as input a list of algorithmic statements which it treats as hypotheses.
It seeks to generate statements which are logically entailed by these hypotheses. 
It outputs a list consisting of the input list together with the newly 
generated statements, if any.
Some rules will require a resource integer $m$ in order to limit the
amount of time that the rule uses.
Resource integers are important in order to allow different rules to take turns being
applied by a supervising process.
Finally, some rules depend on the choice of a library~$\rho$, so
in general a library (or at least an algorithmic sequence which the
rule treats as a library) must be included in the input.

\begin{definition}\label{def42}
A \emph{rule} is an algorithm~$\alpha$ which expects 
an input of the form $[H, \rho, m]$
where $H$ is a list of algorithmic statements, 
$\rho$ is an algorithmic sequence, and $m$ is a positive integer.
For any such input, $\alpha$ is required to halt
with output consisting of a list of algorithmic statements 
containing $H$ as an initial sublist.  Call 
$H$ the \emph{hypothesis list}, $\rho$
the \emph{nominal library}, and
$m$ the \emph{resource integer}. 
The output of the rule $\alpha$ is the \emph{conclusion list}.  

For convenience, we require a \emph{monotonicity property}: 
If $m'\ge m$ and if every item of $H$ is also an item of $H'$,
then every item of the conclusion list for input~$[H, \rho, m]$ is also 
an item of the conclusion list for input~$[H', \rho, m']$.
\end{definition}

\begin{definition}\label{def43}
A \emph{library} is an algorithmic sequence $\rho$ such that $\rho_n$ is a rule for all positive integers $n$.  
\end{definition}

As defined, a library is an infinite sequence of rules, but these rules
are not necessarily distinct.  
In fact, any finite collection of
data can be represented as an algorithmic
sequence $\alpha$ by defining
$\alpha_n$ to be $\alpha_N$ for all $n\ge N$ where $N$ is the
size of the collection. 
Thus the definition does not exclude finite libraries.

\begin{definition}\label{def44}
A rule is \emph{$\rho$-valid} for a library $\rho$ if,
for all hypothesis lists $H$ consisting of only true statements 
and for all resource integers $m$, 
the conclusion list for input $[H, \rho, m]$ consists only of true statements.  
A library $\rho$ is \emph{valid} if it contains only $\rho$-valid rules.  
\end{definition}

\begin{definition}\label{def45}
Let $A_1,\ldots, A_n$ and $B$ be algorithmic statements.  
Let $\rho_k$ be the $k$th rule of a library $\rho$.  
The statement $B$ is a 
\emph{direct $\rho_k$-consequence} of~$A_1,\ldots,A_n$ if there is a
hypothesis list $H$ and a resource integer $m$ such that 
$(i)$ every item of $H$ is in $\{A_1,\ldots,A_n\}$ and 
$(ii)$ $B$ is an item of the conclusion list obtained by 
running $\rho_k$ with input~$[H,\rho,m]$.\footnote{This definition
depends on $\rho$ as well as the particular rule~$\rho_k$.}
\end{definition}

\begin{definition}\label{def46}
A set of algorithmic statements $S$ is \emph{$\rho$-deductively closed} if,
for all $k$, every direct $\rho_k$-consequence of elements in~$S$ 
is itself in~$S$.
\end{definition}

\begin{lemma}\label{lemma47}
If $\rho$ is a valid library then the
set $S$ of true algorithmic statements is $\rho$-deductively closed.
\end{lemma}

\begin{lemma}\label{lemma48}
The intersection of $\rho$-deductive closed sets is 
$\rho$-deductively closed.
\end{lemma}

\begin{definition}\label{def49}
Let $S$ be a set of algorithmic statements.  
The \emph{$\rho$-deductive closure $\overline{S}$ of $S$} is the intersection of all $\rho$-deductively closed sets containing~$S$.
\end{definition}

\begin{lemma}\label{lemma410}
The $\rho$-deductive closure $\overline{S}$ of a set of algorithmic statements is
the minimal $\rho$-deductively closed set containing $S$.  
Thus $\overline{\overline{S}}=\overline{S}$.
\end{lemma}

The $\rho$-deductive closure of a finite set $\{A_1,\ldots,A_n\}$ of algorithmic 
statements  can be explicitly constructed as follows.
Let $f:\bN_{+} \rightarrow \bN_{+} \times \bN_{+}$  be a recursive bijection.  
Define $H_0=[A_1,\ldots,A_n]$.  
For $i>0$, define $H_i$ to be the conclusion list obtained
 by running $\rho_k$ on $[H_{i-1}, \rho, m]$ where $f(i)=(k,m)$. 

\begin{lemma}\label{lemma411}
$B$ is in the $\rho$-deductive closure of $\{A_1,\ldots,A_n\}$ 
if and only if $B$ is an item of $H_i$ for some~$i$.
\end{lemma}

\begin{proof}
Let $S$ be the set of statements on $H_0, H_1, \ldots $. The strategy is to show 
$(i)$ the set of items of any particular $H_i$ is in the $\rho$-deductive closure
(so $S$ is a subset of the $\rho$-deductive closure), and 
$(ii)$ $S$ is $\rho$-deductively closed.

$(i)$ By induction on $i$.  
The case $i=0$ is clear.  
Assume that every item of $H_{i-1}$ is in the $\rho$-deductive closure.  
If $B$ is an item of $H_i$, and if $f(i)=(k,m)$, 
then $B$ is a direct $\rho_k$-consequence of the items of $H_{i-1}$.  
So $B$ is in the $\rho$-deductive closure.

$(ii)$ Suppose $B$ is a direct $\rho_k$-consequence of $C_1,\ldots,C_r \in S$. 
We must show that $B\in S$.
By definition, $B$ is on the conclusion list obtained by running $\rho_k$ with input $[H,\rho,m]$  for some integer $m$ and some list $H$ where
every item of $H$ is in the set $\{C_1,\ldots,C_r\}$.  
By the mononicity requirement for rules,
if $H'$ is any list whose items include each $C_1,\ldots,C_r$ and if $m'\ge m$ 
then $B$ is on the conclusion list when $\rho_k$ is run with input~$[H',\rho,m']$.  

Let $i_0$ be an integer such that $C_1,\ldots, C_r$ are all on $H_{i_0}$.  
There are an infinite number of pairs $(k,m')$ with $m'\geq m$, and all but a finite number are of the form $f(i)$ for $i>i_0$.  
Choose such an $i$.  
So $B$ is on the conclusion list when $\rho_k$ is run with
the input $[H_{i-1},\rho,m']$.  
That is, $B$ is on~$H_i$.  
Thus $B\in S$.
\end{proof}

\begin{definition}\label{def412}
The algorithm $\deduce$ expects an input of the form $[\Gamma, \rho, B]$
where $\Gamma$ is a list of algorithmic statements, 
$\rho$ is a library, and $B$ is an algorithmic statement.  
It computes $H_0,H_1,\ldots$, where $H_0=\Gamma$ and $H_i$ is defined as above.  After computing $H_k$, $\deduce$ checks to see if $B$ is on $H_k$.  
If so, $\deduce$ outputs~$1$;  otherwise, it calculates~$H_{k+1}$.
\end{definition}

\begin{definition}\label{def413}
Let $B$ be an algorithmic statement and $\Gamma$ a list of
algorithmic statements.
The algorithmic statement $\bigl[\deduce, [\Gamma, \rho, B], 1\bigr]$
is denoted as $\Gamma\turnrho B$ 
(usually $\rho$ is a library, but the definition
applies to any datum~$\rho$).
If $\Gamma$ is the list $[A_1,\ldots, A_n]$ one may write 
$A_1, \ldots, A_n\turnrho B$ instead. 
Likewise, $\Gamma, C_1, \dots, C_k \turnrho B$ denotes $\Gamma' \turnrho B$ where $\Gamma'$ is the list obtained by appending $C_1,\ldots, C_k$ to the list~$\Gamma$.  
\end{definition}

\begin{proposition}\label{prop414}
Suppose $\Gamma=[A_1,\ldots,A_n]$ where $A_1,\ldots,A_n$ are algorithmic statements, and suppose $\rho$ is a library.
Then $\Gamma \turnrho B$ if and only if 
$B$ is in the $\rho$-deductive closure of~$\{A_1,\ldots,A_n\}$.
\end{proposition}

\begin{proof}
This follows from Lemma~\ref{lemma411} and the definition of~$\deduce$.
\end{proof}

In particular, if $\Gamma\turnrho B$ then any 
$\rho$-deductively closed set containing all the 
items of~$\Gamma$ contains~$B$.

\begin{corollary} \label{cor415}
Let $A$ and $B$ be algorithmic statements,
$\Gamma$ and $\Gamma'$ lists of algorithmic statements,
and $\rho$ a library.
\begin{list}{}{}
\item[(i)]
If every item of $\Gamma$ is on~$\Gamma'$ and
if $\Gamma \turnrho A$ then $\Gamma' \turnrho A$.
\item[(ii)]
$A \turnrho A$.
\item[(iii)]
If $\Gamma \turnrho A$ and $\Gamma, A \turnrho B$ 
then $\Gamma \turnrho B$.
\item[(iv)]
If $\Gamma\turnrho A$ and if
$\Gamma' \turnrho C_i$ for all items~$C_i$ of~$\Gamma$,
then $\Gamma' \turnrho A$.
\end{list}
\end{corollary} 

\begin{proof}
$(i)$ If $S_1 \subseteq S_2$ then $\overline{S_1} \subseteq \overline{S_2}$.
$(ii)$ $S\subseteq \overline{S}$. $(iii)$ Let $S$ be the $\rho$-deductive
closure of the items of~$\Gamma$. So $A\in S$. Since
$\Gamma, A \turnrho B$ and $S$ is
$\rho$-deductively closed,  $B\in S$. $(iv)$ Let $S$ be the $\rho$-deductive
closure of the items of~$\Gamma'$. So every
item~$C_i$ of~$\Gamma$ is in~$S$. 
Since $\Gamma\turnrho A$ and since
$S$ is deductively closed,  $S$ must contain~$A$.
\end{proof}

\begin{proposition} [Soundness] \label{prop416} 
Suppose every statement on the list~$\Gamma$ is true,  
$\rho$ is a valid library,  and $\Gamma \turnrho B$.
Then $B$ is true.
\end{proposition}

\begin{proof}
The set of true statements $S$ is $\rho$-deductively closed by 
Lemma~\ref{lemma47}. 
The result follows from Proposition~\ref{prop414}.
\end{proof}

Because $\deduce$ is an \emph{internal} predicate
representing deduction, one can use~$\turnrho$ to define
a conditional connective $\imprho$.  (A material
conditional $\dimp$, not dependent on~$\deduce$, will be 
defined in Section~\ref{s12}). 
The algorithm $\deduce$ can also be used
to define an internal provability predicate~$\proverho$.

\begin{definition} \label{def417}
Let $A\imprho B$ denote $A\turnrho B$.
Let $\proverho (A)$ denote~$\truesent\imprho A$.
\end{definition}

The above results, restated in this notation, yield the following.

\begin{proposition} \label{prop418}
Let $A, B$, and $C$ be algorithmic statements, and
$\rho$ a library. Then
\begin{list}{}{}
\item[(i)] $A\imprho A$, and
\item[(ii)]\label{tranprop} if $A\imprho B$ and $B\imprho C$ then $A\imprho C$. 
\end{list}
Moreover, if $\rho$ is a valid library, then
\begin{list}{}{}
\item[(iii)]
if $A\imprho B$ and $A$ are true, then so is~$B$, and
\item[(iv)] if $\proverho(A)$ is true, then so is $A$.
\end{list}
\end{proposition}


\section{Transitivity Rule}\label{s5}

In the next several sections eleven inference rules will be introduced. These rules
will be used to form a stable base (in the sense of Definition~\ref{def131}).
The first is an internal implementation of Proposition~\ref{tranprop}$(ii)$.

\begin{validrule} \label{rule1}
The \emph{Transitivity Rule} is an algorithm that implements the \emph{rule diagram}
$$ \begin{array}{cl}
A \imprho B&\\
B \imprho C&\\
\cline{1-1}
A\imprho C&.
\end{array}
$$
In other words, assuming the input is of the expected form~$[H, \rho, m]$,
the Transitivity Rule first copies the hypothesis list $H$ to a working list $\Delta$. 
Then it looks for a statement of the form $A\imprho B$
and a statement of the form $B \imprho C$ on the hypothesis list~$H$
where $A, B, C$ are algorithmic statements.  
For all such pairs that it finds, the Transitivity Rule appends the 
statement $A \imprho C$ to the working list $\Delta$.  
After processing all such pairs, it outputs the resulting list $\Delta$ as its conclusion list.
\end{validrule}

\begin{proposition}\label{prop51}
The Transitivity Rule is $\rho$-valid for all libraries $\rho$.
\end{proposition}

\begin{proof}
The $\rho$-validity of this rule follows from Proposition~\ref{prop418}$(ii)$.
\end{proof}


\section{Universal Rules}\label{s6}

\begin{validrule}  \label{rule2}
The \emph{Universal Rule} is an algorithm that generates
all true algorithmic statements. 
More specifically,
assuming the input is of the expected form~$[H, \rho, m]$,
the Universal Rule
outputs a list consisting of $H$ appended 
with all $m$-true algorithmic statements.
An algorithmic statement $B$ is \emph{$m$-true} if $(i)$ the datum~$B$ 
has size at most~$m$, $(ii)$ the runtime of the associated process
is at most~$m$, and $(iii)$ $B$ is true.
\end{validrule}

\begin{proposition}\label{prop61}
The Universal Rule is $\rho$-valid for all libraries $\rho$.
\end{proposition}

\begin{proof}
The Universal Rule only appends true statements to the input list.
\end{proof}

\begin{proposition}\label{prop62}
Suppose the library~$\rho$ contains the Universal Rule.
Let $\Gamma$ be a list of algorithmic statements, 
and $A$ and $B$ be algorithmic statements. 
If $B$ is true then $\Gamma \turnrho B$. 
In particular, if $B$ is true, then so is $A \imprho B$ and $\proverho (B)$.
\end{proposition}

\begin{proof}
Every true algorithmic statement is $m$-true for some $m$.
So the $\rho$-deductive closure of any set contains all true statements.
\end{proof}

\begin{corollary}\label{cor63}
If the library $\rho$ is valid and contains the Universal Rule, then an algorithmic statement $A$ is true if and only if $\proverho (A)$.
\end{corollary}

\begin{proof}
This follows from Proposition~\ref{prop418}$(iv)$ and Proposition~\ref{prop62}.
\end{proof}

So, in algorithmic logic, there is a sense in which internal deduction
is complete
for any valid library containing the Universal Rule.  
By Proposition~\ref{prop137}, however,
there is also a sense in which algorithmic logic is inherently incomplete.

\begin{validrule} \label{rule3}
The \emph{Meta-Universal Rule} is an algorithm that implements the rule diagram
$$ \begin{array}{cl}
B\\
\cline{1-1}
A\imprho B&.
\end{array}
$$
More specifically, assuming an input of the expected form~$[H,\rho, m]$, 
the Meta-Universal Rule appends to $H$ all statements of 
the form $A \imprho B$ where 
$(i)$ $B$ is on $H$ and 
$(ii)$ the size of the datum $A \imprho B$ is at most~$m$.
\end{validrule}

The Meta-Universal Rule is the first rule whose validity is contingent on the
\emph{contents} of the library.

\begin{proposition}\label{prop64}
If the library $\rho$ contains the Universal Rule, 
then the Meta-Universal Rule is $\rho$-valid.
\end{proposition}

\begin{proof}
This follows from Proposition~\ref{prop62}.
\end{proof}

\begin{proposition}\label{prop65}
If the library $\rho$ contains the Transitivity Rule and
the Meta-Universal Rule, then 
\begin{list}{}{}
\item[(i)]
\; $A, A \imprho C \; \turnrho \; B \imprho C$, and
\item[(ii)]
\; $A, A \imprho C \; \turnrho \; \proverho (C).$
\end{list}
\end{proposition}

\begin{proof}
$(i)$ Let $S$ be the $\rho$-deductive closure of (the set consisting of) $A$ and $A \imprho C$.  By the Meta-Universal Rule, $B \imprho A$ is in $S$.  By the
Transitivity Rule, $B \imprho C$ is in $S$.

$(ii)$ This is a special case of Part$(i)$ where $B$ is $\truesent$.
\end{proof}


\section{Conjunction} \label{s7}

\begin{definition}\label{def71}
The algorithm $\asand$ expects as input a list $[A,B]$ where $A$ and $B$ are algorithmic statements.  
If $A$ and $B$ are true, then $\asand$ outputs~$1$.  
If either is directly false, then $\asand$ outputs~$0$.  
Otherwise, $\asand$ does not halt.  
If $A$ and $B$ are algorithmic statements,
then $\bigl[ \asand , [A, B], 1 \bigr]$ is denoted
by~$A\wedge B$.

If $\Gamma = [C_1, \ldots , C_k]$ is a list of algorithmic
statements, then the \emph{conjunction}
$C_1\wedge \ldots \wedge C_k$ of $\Gamma$ is 
defined to be $(C_1\wedge \ldots \wedge C_{k-1})\wedge C_k$.
If $k=1$ then the conjunction is simply defined to be $C_1$,
and if $k=0$ (so $\Gamma$ is the empty list) 
then the conjunction is defined to be~$\truesent$.
Observe that the conjunction
$C_1 \wedge \cdots \wedge C_k$ is true if and only if each $C_i$ is true.  
Similarly, the conjunction is directly false if and only if some $C_i$ is directly false. 
\end{definition}

\begin{validrule} \label{rule4}
The \emph{Conjunction Rule} is an algorithm that 
simultaneously implements the following three rule diagrams:
$$
\begin{array}{cccccl}
A \\
B&\quad&A\wedge B&\quad&A\wedge B \\
\cline{1-1}\cline{3-3}\cline{5-5}
A\wedge B && A && B&.
\end{array}
$$
More specifically, 
for any statements $A$ and $B$ on~$H$, 
the Conjunction Rule appends $A\wedge B$ to~$H$.
In addition, for any statement $A\wedge B$ on the given~$H$, 
the Conjunction Rule appends $A$ and $B$ to~$H$.
\end{validrule}

\begin{proposition}\label{prop72}
The Conjunction Rule is $\rho$-valid for all libraries $\rho$.
\end{proposition}

\begin{proposition}\label{prop73}
Let $S$ be a $\rho$-deductively closed set of algorithmic statements 
where $\rho$ is a library containing the Conjunction Rule.  
Let $A_1, \ldots , A_k$ be algorithmic statements where $k\ge 1$.  
The conjunction $A_1 \wedge \cdots \wedge A_k$ is in $S$ if and only 
if each $A_i$ is in~$S$. 
(If $k=0$ assume that $\rho$ contains the Universal Rule instead of the
Conjunction Rule).
\end{proposition}

\begin{corollary}\label{cor74}
Let $\rho$ be a library containing the Conjunction Rule.
The logical connective $\wedge$ satisfies both the symmetry and associativity laws:
\begin{list}{}{}
\item[(i)]
\; $A \wedge B \; \turnrho \; B \wedge A$.
\item[(ii)]
\; $(A \wedge B) \; \wedge C \; \turnrho \; A \wedge (B \wedge C)$
\quad and \quad
$A \wedge (B \wedge C)\; \turnrho \; (A \wedge B) \wedge C $.
\end{list}
\end{corollary}

\begin{corollary}\label{cor75}
Let $\Gamma = [C_1, \ldots , C_k]$ be a list of algorithmic
statements, 
$A_1, \ldots , A_n, B$ be algorithmic statements,
and $C = C_1 \wedge \cdots \wedge C_k$.
Assume that $\rho$ contains the Conjunction Rule
and the Universal Rule (for the case $k=0$ or $n=0$).  
Then 
\begin{list}{}{}
\item[(i)]
$\Gamma \; \turnrho \; A_1 \wedge \cdots \wedge A_n$
\quad
if and only if
\quad
$\Gamma \; \turnrho \; A_i$ for each $A_i$, and 
\item[(ii)]
$\Gamma \turnrho B$ \quad if and only if \quad $C \imprho B$.
\end{list}                                              
\end{corollary}

\begin{validrule}  \label{rule5}
The \emph{Meta-Conjunction Rule} is an algorithm that implements the rule diagram
$$ \begin{array}{cl}
A \imprho B\\
A \imprho C\\
\cline{1-1}
A\imprho (B \wedge C)&.
\end{array} $$
\end{validrule}

\begin{proposition}\label{prop76}
The Meta-Conjunction Rule is $\rho$-valid for all libraries $\rho$ containing the Conjunction Rule.
\end{proposition}

\begin{proof}
Assume $A \imprho B$ and $A \imprho C$.  Let $S$ be 
the $\rho$-deductive closure of $A$.  By assumption $B$ and $C$
are in $S$.  By the Conjunction Rule $B \wedge C$ is also in $S$.
Therefore, $A \imprho B \wedge C$.
\end{proof}

Several laws can be deduced from the above rules.

\begin{proposition}\label{prop77}
If $\rho$ contains all the above rules, then
\begin{list}{}{}
\item[(i)]
\; $A \imprho B \;\; \turnrho \; \; A \imprho B \wedge A$,
\item[(ii)]
\; $A \imprho B, \; B \wedge A  \imprho C \; \; \turnrho \;\; A \imprho C$, and
\item[(iii)]
\; $A \imprho B \;\; \turnrho \;\; C \wedge A \imprho C \wedge B$,
\quad
$A \imprho B \;\; \turnrho \;\; A \wedge C \imprho B \wedge C$.
\end{list}
\end{proposition}

\begin{proof}
$(i)$  Let $S$ be the $\rho$-deductive closure of $A \imprho B$.  
The statement $A \imprho A$ is true by Proposition~\ref{prop418}$(i)$.  
By the Universal Rule, $A \imprho A$ is in $S$.  
So by the Meta-Conjunction Rule $A \imprho B \wedge A$ is in $S$. 

$(ii)$  Let $S$ be the $\rho$-deductive closure 
of~$ A \imprho B$ and $B \wedge A \imprho C$.  
By the first part, $A \imprho B \wedge A$ is in $S$.  
So, by the Transitivity Rule, $A \imprho C$ is in~$S$. 

$(iii)$ This follows by a similar argument.
\end{proof}

\begin{proposition}\label{prop78}
Suppose $\rho$ contains all the above rules. 
If $B \wedge A \imprho C$ then $B \imprho (A \imprho C)$.
\end{proposition}

\begin{proof}
Let $S$ be the $\rho$-deductive closure of $B$.  
By supposition $B \wedge A \imprho C$ holds, so
is in~$S$ by the Universal Rule.  
By the Meta-Universal Rule, $A \imprho B$ is in~$S$.  
Finally, by Proposition~\ref{prop77}$(ii)$, $A \imprho C$ is in~$S$.
\end{proof}

\begin{theorem}\label{theorem79}
Suppose $\rho$ contains all the above rules.  
Let $\Gamma$ be a list of algorithmic statements, 
and let $A$ and $C$ be algorithmic statements. 
If $\Gamma, A \; \turnrho\; C$ then $\Gamma \; \turnrho \; A\imprho C$.
\end{theorem}

\begin{proof}
Let $\Gamma=[B_1, \ldots, B_k]$ and 
$B=B_1 \wedge \cdots \wedge B_k$.
If $\Gamma, A \turnrho C$, then
$B\wedge A \imprho C$ by Corollary~\ref{cor75}$(ii)$.
By Proposition~\ref{prop78}, $B\imprho (A\imprho C)$ holds.
By Corollary~\ref{cor75}$(ii)$ again, $\Gamma \turnrho A\imprho C$.
\end{proof}

\begin{corollary}\label{cor710}
If $\rho$ contains all the above rules, then
\begin{list}{}{}
\item[(i)]
\; $A \;\; \turnrho \;\; B \imprho A \wedge B$,
\item[(ii)]
\; $A \imprho B \; \;  \turnrho \; \; (B \imprho C) \imprho (A \imprho C)$,
\item[(iii)]
\; $A \imprho B \; \;  \turnrho \; \; (C \imprho A) \imprho (C \imprho B)$, and
\item[(iv)]
\; $B \wedge A \imprho C\; \; \turnrho \;\;  B \imprho (A \imprho C) $.
\end{list}
\end{corollary}

\begin{proof}
$(i)$ By the Conjunction Rule, $A, B \; \turnrho \; A \wedge B$. 
Now use Theorem~\ref{theorem79}.  

$(ii)$ By the Transitivity Rule, $A \imprho B, B\imprho C \; \turnrho \; A \imprho C$.  
Now use Theorem~\ref{theorem79}. Part$(iii)$ is similar.

$(iv)$ Let $S$ be the $\rho$-deductive closure of $B \wedge A \imprho C$ and $B$.  
By the Meta-Universal Rule, $A \imprho B$ is in $S$.  
By Proposition~\ref{prop77}$(ii)$, $A \imprho C$ is in~$S$.  
Thus $B \wedge A \imprho C, B \; \turnrho \; A \imprho C$.  
Now use Theorem~\ref{theorem79}.
\end{proof}


\section{Biconditional}\label{s8}

\begin{definition}\label{def81}
Let $A \iffrho B$ denote $(A \imprho B) \wedge (B \imprho A)$.
\end{definition}

\begin{proposition}\label{prop82}
The following laws hold for any library $\rho$:
\begin{list}{}{}
\item[(i)]
\; $A \iffrho A$,
\item[(ii)]
\; If $A \iffrho B$ then $B \iffrho A$, and
\item[(iii)]
\; If $A \iffrho B$ and $B \iffrho C$, then $A \iffrho C$.
\end{list}
\end{proposition}

Some results concerning conjunction can be conveniently expressed
with the biconditional.

\begin{proposition}\label{prop83}
Suppose $\rho$ is a library containing the Conjunction and Universal Rules.
Then
\begin{list}{}{}
\item[(i)]
\; $A \iffrho A \wedge A$,
\item[(ii)]
\; $A \iffrho A \wedge \truesent$,
\item[(iii)]
\; $A \wedge B \iffrho B \wedge A$, and
\item[(iv)]
\; $A \wedge (B \wedge C) \iffrho (A \wedge B) \wedge C$.
\end{list}
\end{proposition}


\section{Disjunction}\label{s9}

\begin{definition}\label{def91}
The algorithm $\asor$ expects as input a list $[A,B]$ where $A$ and $B$ are algorithmic statements.  
If either $A$ or $B$ are true, then $\asor$ outputs~$1$.  
If both are directly false, then $\asor$ outputs~$0$. 
Otherwise, $\asor$ does not halt.

If $A$ and $B$ are algorithmic statements, 
then we denote $\bigl[\asor, [A, B], 1\bigr]$ by $A \vee B$. 
The statement $A \vee B$ is true if and only if either $A$ is true or $B$ is true.  
Similarly, $A \vee B$  is directly false if and only if both $A$ and $B$ are directly false.
\end{definition}

\begin{validrule} \label{rule6}
The \emph{Disjunction Introduction Rule} is an algorithm that simultaneously implements the following two rule diagrams:
$$
\begin{array}{cccl}
A&&B\phantom{.}\\
\cline{1-1}\cline{3-3}
A\vee B &\quad\quad& A \vee B&.
\end{array}
$$
More specifically,
assuming an input in the expected form~$[H,\rho, m]$, 
the Disjunction Introduction Rule appends to $H$ all statements of the 
form $A \vee B$ where 
(i) either $A$ or $B$ is on $H$ and 
(ii) the size of $A \vee B$ is at most~$m$.
\end{validrule}

\begin{proposition}\label{prop92}
The Disjunction Introduction Rule is $\rho$-valid for all libraries~$\rho$.
\end{proposition}

At this point one might expect an algorithmic
disjunction elimination rule allowing the deduction of
$C$ from $A \imprho C$, $B \imprho C$, and $A \vee B$.
The difficulties of such a rule will be 
discussed in Section~\ref{s14}.
An unproblematic but weaker version of this rule
can be produced by requiring a sort of verification
for the hypotheses $A \imprho C$ and $B \imprho C$.
The following rule implements this idea.\footnote{To see
the usual disjunction elimination rule, 
think of $G$ as $\truesent$.
Allowing general~$G$ is important in the proof of Proposition~\ref{prop96}.}

\begin{validrule}  \label{rule7}
The \emph{Disjunction Elimination Rule} is an algorithm, 
denoted $\disjunctelim$, that implements the rule diagram
$$
\begin{array}{lc}
G \wedge A \imprho C & *\\
G \wedge B \imprho C & *\\
G &\\
A \vee B &\\
\cline{1-1}
C&
\end{array}
$$
where $*$ indicates that the corresponding statement
must be \emph{verified}.  
More specifically, 
assuming an input of the expected form~$[H,\rho, m]$, 
whenever $\disjunctelim$ finds four
statements on~$H$ of the form of the premises of the rule diagram,
it determines if the runtimes of the processes
associated with the first two statements in the diagram are 
less than~$m$. 
If the runtimes are both less than~$m$ and if
both statement are \emph{true},
then $\disjunctelim$ appends the statement represented
by $C$ to the conclusion list.
\end{validrule}

Proposition~\ref{prop418}\textit{(iii)} and the 
algorithmic definition of~$\vee$ gives validity:

\begin{proposition}\label{prop93}
The Disjunction Elimination Rule is $\rho$-valid for all valid
libraries $\rho$.
\end{proposition}

This is the first rule we have considered where the validity of the rule is
contingent on the validity of the library.  If
the library $\rho$ is valid, then this rule is $\rho$-valid,
but in Section~\ref{s14} we shall see several examples of valid rules
that cannot themselves be contained in a stable base
(in the sense of Definition~\ref{def131}).
Since the goal is to form a stable base of inference rules,
we need something stronger than
the above proposition. Theorem~\ref{theorem95}
is sufficient.

\begin{lemma}\label{lemma94}
If a library~$\rho$ is not valid, but contains the Conjunction Rule,
then there are algorithmic statements $A$ and $B$
such that $A$ and $A\imprho B$ are true, but $B$ is false.
\end{lemma}

\begin{proof}
Since $\rho$ is not valid, there is a rule $\rho_k$ in~$\rho$
which is not $\rho$-valid.  In other words, there is 
a list $H$ of true statements and an integer $m$ 
such that when the list $[H,\rho,m]$ is given as input to $\rho_k$, 
the rule generates an output list containing at least one false statement $B$.  
Let $H=[A_1, \ldots, A_n]$ and let $A = A_1 \wedge \cdots \wedge A_n$.  
Note that $A$ is true.
Let $S$ be the $\rho$-deductive closure of~$A$. 
Since $\rho$ contains the Conjunction Rule, each $A_i$ is in~$S$.
So $B$ is in~$S$ by the definition of the deductive closure.  
Thus $A \imprho B$ is true.
\end{proof}

\begin{theorem}\label{theorem95}
Let $\rho$ be a library that contains at least the Conjunction Rule 
and the Disjunction Elimination Rule.  
Suppose that every rule in $\rho$ other than the Disjunction
Elimination Rule is $\rho$-valid.  
Then $\rho$ is valid.
\end{theorem}

\begin{proof}
Suppose to the contrary that $\rho$ is not valid.  
By the previous lemma there are statements $D$ and $E$ such that $D$ is true, 
$D \imprho E$ is true, but $E$ is false.  
Choose $D$ and $E$ so that the runtime $r$ of the process associated
with $D \imprho E$ is minimal.

Let $H_0 = [D]$.  
Since $D \imprho E$, when $[H_0, \rho, E]$ is input to $\deduce$
the output is~$1$.  Recall that $\deduce$ generates a monotonic sequence $H_0, H_1, \ldots $ of lists, 
and since it outputs~$1$, it eventually generates a list $H_k$ containing~$E$.  
Thus, since $H_0$ contains only true statements but $E$ is false,
there is a unique $i \ge 1$ such that $H_{i-1}$ contains only 
true statements and $H_i$ contains at least  one false statement~$C$.  
Let the function $f$ be as in the definition of $\deduce$, and let $f(i)= (k,m)$.  
Thus $H_i$ is obtained by running $\rho_k$ with input $[H_{i-1}, \rho, m]$.  
Note that $\rho_k$ cannot be $\rho$-valid, so $\rho_k$ 
must be $\disjunctelim$ (since we assumed that all other rules are
$\rho$-valid).  
Since $\disjunctelim$ generates~$C$, 
$H_{i-1}$ must contain statements of the form 
$(i)$~$G \wedge A \imprho C$, 
$(ii)$~$G \wedge B \imprho C$, 
$(iii)$~$G$, and 
$(iv)$~$A \vee B$.  
These four statements are true since they are on $H_{i-1}$.  
So either $A$ or $B$ is true, and it is enough to consider the case
where $A$ is true. In this case $G \wedge A$ is true.  
Since $\disjunctelim$ generates the statement $C$, 
it must first run the process associated with $G \wedge A \imprho C$
and determine that the statement is true.  
The runtime $r$ of the global process associated with $D \imprho E$ 
must be strictly larger than the runtime~$r'$ associated 
\hbox{with~$G \wedge A \imprho C$}
(since $r'$ is the runtime a subprocess of a subprocess of the global process 
associated with~$D \imprho E$).  
Since $r' < r$ and since both $G \wedge A$ and $G \wedge A \imprho C$ are true, 
it follows from the definition of $r$ that $C$ must be true, a contradiction.
\end{proof}

\begin{proposition}\label{prop96}
Let $\rho$ be a library containing the Universal, Conjunction, 
and Disjunction Elimination Rules.
\begin{list}{}{}
\item[(i)]
\; If $G \wedge A \imprho C$ and $G \wedge B \imprho C$ then
$G \wedge (A \vee B) \imprho C$.
\item[(ii)]
\; If $\;\Gamma, A\; \turnrho\; C\;$ and $\;\Gamma, B \;\turnrho\; C$,
 then $\;\Gamma, A \vee B \; \turnrho\; C$.
\end{list}
\end{proposition}

\begin{proof}
$(i)$ Let $S$ be the $\rho$-deductive closure of $G \wedge (A \vee  B)$.  
We must show that $C$ is in $S$.  
By the Conjunction Rule, $G$ and $A \vee B$ are in $S$.
By the Universal Rule, $G \wedge A \imprho C$ and $G \wedge B \imprho C$ 
are also in~$S$.  
So by the Disjunction Elimination Rule, $C$ is in~$S$ 
(where $\disjunctelim$ needs a resource number~$m$ larger than the runtimes 
associated with~$G \wedge A \imprho C$ and $G \wedge B \imprho C$).

 $(ii)$ This follows from Part$(i)$ and Corollary~\ref{cor75}$(ii)$.
\end{proof}

\begin{proposition}\label{prop97}
If $\rho$ contains all the above rules, then
\begin{list}{}{}
\item[(i)]
\; $A \vee \truesent \iffrho  \truesent$,
\item[(ii)]
\; $A \iffrho A \vee A$,
\item[(iii)]
\; $A \vee B \iffrho B \vee A$, and
\item[(iv)]
\; $A\iffrho A \wedge (A \vee B)$.
\item[(v)]
\; $A \iffrho A \vee (A \wedge B)$.
\item[(vi)]
\; $A \vee (B \vee C) \iffrho (A \vee B) \vee C$.
\end{list}
\end{proposition}

\begin{proof}
(\emph{i}) to (\emph{v}) These are similar to and easier than Part~(\emph{vi}):

(\emph{vi}) The Disjunction Introduction Rule (twice)
gives \hbox{$B \turnrho \; (A \vee B) \vee C$}.  Likewise,
$C \; \turnrho \; (A \vee B) \vee C$. Proposition~\ref{prop96}$(ii)$
gives $B \vee C \; \turnrho \; (A \vee B) \vee C$.  The
Disjunction Introduction Rule (twice) 
gives $A \turnrho \; (A \vee B) \vee C$. 
Finally, Proposition~\ref{prop96}$(ii)$
gives $A \vee (B \vee C) \; \turnrho \; (A \vee B) \vee C$. 
This gives one direction.
The other direction follows from a similar argument.
\end{proof}

\begin{proposition}\label{prop98}
If $\rho$ contains all the above rules, then
\begin{list}{}{}
\item[(i)]
\; $A \wedge (B \vee C) \iffrho (A \wedge B) \vee (A \wedge C)$, and
\item[(ii)]
\; $A \vee (B \wedge C) \iffrho ( A \vee B) \wedge (A \vee C)$.
\end{list}
\end{proposition}

\begin{proof}
$(i)$  By the Disjunction Introduction Rule,
$A \wedge B \; \turnrho \; (A \wedge B) \vee (A \wedge C)$ and
$A \wedge C \; \turnrho \; (A \wedge B) \vee (A \wedge C)$.  Now use
Proposition~\ref{prop96}$(i)$ to show
$A \wedge (B \vee C) \turnrho (A \wedge B) \vee (A \wedge C)$.
The other direction is similar.

$(ii)$ Showing
$A \vee (B \wedge C) \turnrho ( A \vee B) \wedge (A \vee C)$
is similar to Part$(i)$.
For the other direction, first show
$C, A \turnrho A \vee (B \wedge C)$ and $C, B \turnrho A \vee (B \wedge C)$ using the
Disjunction Introduction and Conjunction Rules.  
Use Proposition~\ref{prop96}$(ii)$ to get
$C, A \vee B \; \turnrho \; A \vee (B \wedge C)$.  In other words,
$A \vee B, C\; \turnrho \; A \vee (B \wedge C)$.  
Use the Disjunction Introduction Rule to get
$A \vee B, A \; \turnrho \; A \vee (B \wedge C)$.  
Use Proposition~\ref{prop96}$(ii)$ again to get
$A \vee B, A\vee C\; \turnrho \; A \vee (B \wedge C)$.  
Finally, use
Proposition~\ref{cor75}$(ii)$ to get the conclusion.
\end{proof}

\begin{validrule} \label{rule8}
The \emph {Meta-Disjunction Rule} is an algorithm that implements the rule diagram
$$ \begin{array}{ll}
G \wedge A \imprho C\\
G \wedge B \imprho C\\
\cline{1-1}
G \wedge (A \vee B) \imprho C&.\\
\end{array}
$$
\end{validrule}

\begin{proposition}\label{prop99}
If the library $\rho$ contains the
Universal, Conjunction, and Disjunction Elimination Rules,
then the Meta-Disjunction Rule is $\rho$-valid.
\end{proposition}

\begin{proof}
This follows from Proposition~\ref{prop96}$(i)$.
\end{proof}

\begin{proposition}\label{prop910}
If $\rho$ contains all the above rules, then
$$
A \imprho C, B \imprho C \; \turnrho \; A \vee B \imprho C.
$$
\end{proposition}

\begin{proof}
Let $S$ be the $\rho$-deductive closure of the two hypotheses.
By the Conjunction, Universal, and Transitivity Rules,
$\truesent \wedge A \imprho C$ and $\truesent \wedge B \imprho C$
are in~$S$.
By the Meta-Disjunction Rule, $\truesent \wedge (A \vee B) \imprho C$
is in $S$.
By the Universal Rule, $\truesent$ is in $S$.  So, by
Corollary~\ref{cor710}$(i)$, $A \vee B \imprho \truesent \wedge (A \vee B)$
is in $S$.
Finally, by the Transitivity Rule, $A \vee B \; \imprho \; C$ is in $S$.
\end{proof}

\begin{proposition}\label{prop911}
Assume that $\rho$ contains all of the rules defined above.  Then
$A \imprho B \;\; \turnrho \;\;  C \vee A \imprho C \vee B$ and
$A \imprho B \;\; \turnrho \; \;  A \vee C  \imprho B \vee C$.
\end{proposition}

\begin{proof}
Let $S$ be the deductive closure of $A\imprho B$.
By the Disjunction Introduction, Universal, and the Transitivity
Rules,  $A\imprho C\vee B$ and $C\imprho C\vee B$
are in~$S$. By Proposition~\ref{prop910},
$C \vee A \imprho C \vee B$ is in~$S$.
Similarly,  $A \vee C \imprho B \vee C$ is in~$S$.
\end{proof}


\section{Negation}\label{s10}

\begin{definition}\label{def101}
Let $A$ be an algorithmic statement.  
The statement $\notrho A$ is defined to be $A \imprho \falsesent$.
\end{definition}

\begin{proposition}\label{prop102} 
If $\rho$ contains all of the above rules, then
\begin{list}{}{}
\item[(i)]
\;\;  $A \imprho B, \; \notrho B \;\; \turnrho \;\; \notrho A$, 
\item[(ii)]
\;\;  $A \imprho B \;\; \turnrho \;\; \notrho B \imprho \; \notrho A$, and
\item[(iii)]
\;\;   $A, \; \notrho A \;\; \turnrho \;\; \notrho B$.
\end{list}
\end{proposition}

\begin{proof}
$(i)$ Use the Transitivity Rule.  

\noindent
$(ii)$ Use Part$(i)$ and Theorem~\ref{theorem79}.  

\noindent
$(iii)$ Use the Meta-Universal Rule to form 
$B \imprho A$. Then use Part$(i)$.
\end{proof}

One might expect the law $A, \notrho A \; \turnrho \; B$ to hold.
Unfortunately it often fails.
The instability of the corresponding rule will be
discussed in Section~\ref{s14}.

\begin{proposition}[De Morgan]\label{prop103} 
If $\rho$ contains all of the above rules, then 
\begin{list}{}{}
\item[(i)]
\;\; $\notrho (A \vee B) \iffrho \; \notrho A \, \wedge \notrho B$, and
\item[(ii)]
\;\; $\notrho A \vee \notrho B \;\; \turnrho \; \notrho (A \wedge B)$.
\end{list}
\end{proposition}

\begin{proof}
$(i)$  Let $S$ be the $\rho$-deductive closure of $\notrho (A \vee B)$.
In other words, $A \vee B \imprho \falsesent$ is in $S$.  Use the
Disjunction Introduction Rule to get $A \imprho A \vee B$ and the
Universal Rule to show that it is in $S$.
So, $A \imprho \falsesent$ is in $S$ by the Transitivity Rule.  In other words, 
$\notrho A$ is in $S$.  Likewise, $\notrho B$ is in $S$.  The Conjunction
Rule gives that $\notrho A  \, \wedge \notrho B$ is in $S$. So
 $\notrho (A \vee B) \; \turnrho \; \notrho A \, \wedge \notrho B$.

For the other direction, let $S$ be the deductive closure of
$\notrho A \, \wedge \notrho B$.  By the Conjunction Rule
$A \imprho \falsesent$ and
$B \imprho \falsesent$ are in $S$.  By Proposition~\ref{prop910},
$A \vee B \imprho \falsesent$ is in $S$.  So
$\notrho  A \, \wedge \notrho B \; \turnrho \notrho (A \vee B)$.

$(ii)$ Let $S$ be the $\rho$-deductive closure of $A \imprho \falsesent$.
Use the Conjunction Rule to get $A \wedge B \imprho A$ and the
Universal Rule to show that it is in $S$.  So, by the Transitivity Rule,
$A \wedge B \imprho \falsesent$ is in $S$.

 Therefore, $\notrho A \; \turnrho \notrho (A \wedge B)$.  Similarly,
$\notrho B \; \turnrho  \notrho (A \wedge B)$.  So, by
Proposition~\ref{prop96}$(ii)$, 
$\notrho \, A \vee \notrho B \; \turnrho \; \notrho (A \wedge B)$
\end{proof}

The problems with the full  converse
$\notrho (A \wedge B) \imprho \; \notrho A  \vee \notrho B$
of the second part of De Morgan will be addressed 
in Theorem~\ref{theorem1412}.  
Part$(ii)$ of the following gives a partial version.

\begin{proposition}\label{prop104} 
If $\rho$ contains all the above rules, then
\begin{list}{}{}
\item[(i)]
$\; \notrho (A \wedge B), \, B \; \turnrho \; \notrho A$, and
\item[(ii)]
$\; \notrho (A \wedge B), \, B \vee \notrho B \; \turnrho \; \notrho A \vee \notrho B$. 
\end{list}
\end{proposition}

\begin{proof}
$(i)$ Let $S$ be the $\rho$-deductive closure of  $A \wedge B \imprho \falsesent$
and $B$.  By Corollary~\ref{cor710}$(i)$, $A \; \imprho \; B \wedge A$ is in $S$.
By Corollary~\ref{cor74}$(i)$, $B \wedge A \; \imprho \; A \wedge B$ holds
so is in $S$ by the Universal Rule.  By applying the Transitivity Rule twice,
$A \imprho \falsesent$ is in $S$.

$(ii)$  
Both
$\notrho (A \wedge B), \, B  \; \turnrho  \notrho A \vee \notrho B$ and
$\notrho (A \wedge B), \, \notrho B \, \turnrho  \notrho A \vee \notrho B$
hold.
The first follows by Part$(i)$ and the Disjunction Introduction Rule.
The second is a consequence of the  Disjunction Introduction Rule.  
So
by Proposition~\ref{prop96}$(ii)$ the conclusion holds.
\end{proof}

One might expect the law $A \vee B, \notrho B \; \turnrho \; A$ to hold.  
Problems with this law will be discussed in
Section~\ref{s14}.  A partial version is given by the following.

\begin{proposition}\label{prop105} 
If $\rho$ contains the above rules, then
$$
\notrho A \vee B, \notrho B \; \turnrho \; \notrho A.
$$
\end{proposition}

\begin{proof}
Corollary~\ref{cor415} gives $\notrho B, \notrho A \; \turnrho \notrho A$.
Proposition~\ref{prop102}$(iii)$ gives 
\hbox{$\notrho B, B \; \turnrho \notrho A$}.
Finally, Proposition~\ref{prop96}$(ii)$ gives
$\notrho B, \notrho A \vee B \; \turnrho \; \notrho A$.
\end{proof}

The proofs of the propositions above are not contingent on any special properties
of the statement $\falsesent$ itself:
similar results can be derived if $\notrho U$ is systematically
replaced with $U\imprho E$ where $E$ is \emph{any}
fixed statement.
The following proposition, however,  uses a property specific to $\falsesent$: 
if $\rho$ contains the Elimination of Case Rule defined below, 
then $\falsesent \imprho B$ holds for any $B$.  

\begin{proposition}\label{prop106} 
Assume that $\falsesent \imprho B$ holds for any $B$ and
that   $\rho$ contains all the above rules.
\begin{list}{}{}
\item[(i)]
\; If $\notrho A$ then $A\imprho B$.
\item[(ii)]
\; $\notrho A \; \turnrho \; A \imprho B$.
\item[(iii)]
\; $A, \notrho A \; \turnrho \; B \imprho C$.
\item[(iv)]
\; If $\notrho A$ then $A \vee B \; \turnrho \; B$.
\item[(v)]
\; $\notrho A \; \turnrho \; A \vee B \imprho B$.
\item[(vi)]
\; $\notrho A \vee B \; \turnrho \; A \imprho B$.
\end{list}
\end{proposition}

\begin{proof}
$(i)$ By assumption, $A \imprho \falsesent$ and $\falsesent \imprho B$.
So, the conclusion follows from Proposition~\ref{prop418}$(ii)$.

$(ii)$  Let $S$ be the $\rho$-deductive closure of $A \imprho \falsesent$.
Since $\falsesent \imprho B$ holds, it is in $S$ by the Universal Rule.
So, by the Transitivity Rule, $A \imprho B$ is in $S$.

$(iii)$  Let $S$ be the $\rho$-deductive closure of $A$ and $\notrho A$.
Use Part$(ii)$ to get $A \imprho C$ in $S$.  By the Meta-Universal
Rule, $B \imprho A$ is in $S$.  So, by
the Transitivity Rule, $B \imprho C$ is in $S$.

$(iv)$  Assume $\notrho A$.  So by Part$(i)$, $A \turnrho B$.  Since
$B \turnrho B$, the conclusion follows from 
Proposition~\ref{prop96}$(ii)$.  

$(v)$ Let $S$ be the $\rho$-deductive closure of $\notrho A$.  By
Part$(ii)$, $A \imprho B$ is in $S$.  By the Universal Rule, $B \imprho B$
is in $S$.  By Proposition~\ref{prop910}, $A \vee B \imprho B$ is in $S$.

$(vi)$ By Part$(ii)$ above, $\notrho A \; \turnrho \; A \imprho B$.  By the
Meta-Universal Rule, $B \; \turnrho \; A \imprho B$.  The
conclusion follows from Proposition~\ref{prop96}$(ii)$.
\end{proof}

Part$(iii)$ above is a weak version of the ideal law $A, \notrho A \; \turnrho \; B$.
Part$(iv)$ and Part$(v)$ are closely related to the ideal law $A \vee B, \notrho A \;
 \turnrho \; B$.  
And the converse $A \imprho B \; \turnrho \notrho A \vee B$ of
Part$(vi)$ is another ideal law.  
The instability of the corresponding rules
is addressed in Section~\ref{s14}.


\section{Strong Negation}\label{s11}

\begin{definition}\label{def111}
The algorithm $\sneg$ expects as input an algorithmic 
statement $[\alpha, u, v]$.  
It runs $\alpha$ as a subprocess with input $u$.  
If this subprocess halts with output $v$, then $\sneg$ outputs $0$.  
If the subprocess halts with output other than $v$, then $\sneg$ outputs $1$.  
Otherwise $\sneg$ does not halt.

Let $A$ be an algorithmic statement.  
Then the \emph{strong negation} of $A$, denoted by~$-A$,
is the algorithmic statement $[\sneg, A, 1]$.
Note that $-A$ is true if and only if $A$ is directly false, 
and $-A$ is directly false if and only if $A$ is true.  
In particular, $-\falsesent$ is true and $-\truesent$ is directly false.
\end{definition}

\begin{validrule} \label{rule9}
The \emph{Elimination of Case Rule} is an algorithm that
implements the rule diagram
$$ \begin{array}{cl}
A\vee B\\
-A\phantom{-}\\
\cline{1-1}
\phantom{-}B\phantom{-}&.\\
\end{array}
$$
\end{validrule}

\begin{proposition}\label{prop112}
The above rule is $\rho$-valid for all libraries~$\rho$.
\end{proposition}

\begin{proposition}\label{prop113}
If $\rho$ contains all the above rules then
\begin{list}{}{}
\item[(i)]
\; $A, -A \; \turnrho \; B$,
\item[(ii)]
\; $\falsesent \; \turnrho \; B$,
\item[(iii)]
\; $-A \; \turnrho \notrho A$, and
\item[(iv)]
\; $\falsesent \vee A \iffrho A$ \quad and
 \quad  $\falsesent \wedge A \iffrho \falsesent$.
\end{list}
\end{proposition}

\begin{proof}
$(i)$ Let $S$ be the $\rho$-deductive closure of $A$ and $-A$.
By the Disjunction Introduction Rule, $A \vee B$ is in $S$.  So,
by the Elimination of Case Rule, $B$ is in $S$.  

$(ii)$  By the Universal Rule, $\falsesent \turnrho -\falsesent$.  By
Part$(i)$, $\falsesent, - \falsesent \turnrho B$.  The conclusion
follows by Corollary~\ref{cor415}$(iii)$.
  
$(iii)$ By Part$(i)$, $-A, A \; \turnrho \; \falsesent$.  By Theorem~\ref{theorem79},
$-A \; \turnrho \; A \imprho \falsesent$.

 $(iv)$ The first biconditional follows
 from the Disjunction Introduction Rule, Part$(ii)$, and
 Proposition~\ref{prop96}$(ii)$.  The second follows from the Conjunction Rule,
 Part$(ii)$, and Proposition~\ref{cor75}$(i)$.
 \end{proof}

\begin{validrule} \label{rule10}
The \emph{Double Negation Rule} is an algorithm that 
simultaneously implements the following rule diagrams:
$$ \begin{array}{cccl}
A&&--A\\
\cline{1-1}\cline{3-3}
--A&\quad&A&.\\
\end{array}
$$
\end{validrule}

\begin{validrule} \label{rule11}
The \emph{Strong De Morgan Rule} is an algorithm that
simultaneously implements the following
rule diagrams:
$$ \begin{array}{cccccccl}
-(A\vee B)&& -A \wedge -B &&-(A\wedge B)&& -A\vee-B\\
\cline{1-1}\cline{3-3}\cline{5-5}\cline{7-7}
-A\wedge-B&&-(A\vee B)&&-A\vee-B&&-(A\wedge B)&.\\
\end{array}
$$
\end{validrule}

\begin{proposition}\label{prop114}
The Double Negation and Strong De Morgan Rules
are $\rho$-valid for all libraries~$\rho$.
\end{proposition}

\begin{proposition}\label{prop115}
If $\rho$ contains the Double Negation and Strong De Morgan Rules,
then
\begin{list}{}{}
\item[(i)]
\; $A \iffrho --A$,
\item[(ii)]
\; $- (A \vee B) \iffrho -A \wedge -B$, \quad and \quad
$- (A \wedge B) \iffrho -A \vee -B$.
\end{list}
\end{proposition}

 
\section{Material Conditional}\label{s12}

\begin{definition}\label{def121}
Define the \emph{material conditional} $A\dimp B$ to be $-A \vee B$.
Define $\ahalt(A)$ to be $-A \vee A$. 
\end{definition}

Note that $\ahalt(A)$ is $A \dimp A$.  Also note that the statement $\ahalt(A)$ is true
if and only if the process associated with $A$ halts.

\begin{proposition}\label{prop122}
If $\rho$ contains all the above rules, then $ A \dimp B \; \turnrho \; A \imprho B$.
\end{proposition}

\begin{proof}
Use Proposition~\ref{prop113}$(i)$ and Theorem~\ref{theorem79} to get
$-A \turnrho A \imprho B$.  
By the Meta-Universal rule,
$B \turnrho A \imprho B$.  
Finally, use Proposition~\ref{prop96}$(ii)$. 
\end{proof}

The material conditional $\dimp$ has many of the properties
one would expect.  Indeed, some of its properties are stronger than
those of the connective~$\imprho$.  
The connective~$\dimp$ 
has, however, several striking weaknesses.  
The following are true for the material conditional
in classical logic: 
$$A \dimp A, \; 
A \dimp A \vee B, \;
A \wedge B \dimp A, \;
A \dimp (B \dimp A), \;
(A \dimp B) \wedge (B \dimp C) \dimp (A \dimp C),
$$
$$
A \wedge (A \dimp B) \dimp B, \;
\text{and} \;\;
(A \vee B) \wedge (A \dimp C) \wedge (B \dimp C) \dimp C.
$$ 

\smallskip
\noindent
But if $A$, $B$, and $C$ are chosen so that
$\ahalt(A)$, $\ahalt(B)$ and $\ahalt(C)$ are false, 
then these statements are all false for the material conditional
of Definition~\ref{def121}.

Some such tautologies of classical logic can, however, be interpreted
to form corresponding laws of algorithmic logic containing
the connective~$\imprho$ (or equivalently~$\turnrho$) or containing
a mixture of both~$\dimp$ and~$\imprho$
(or~$\turnrho$).  For example, $A \imprho A$, \; $A \imprho A \vee B$, 
and $A \wedge B \imprho A$ hold in general for libraries $\rho$ containing
all the above rules. 
The following give further examples (Parts$(i)$ and $(iv)$-$(vi)$ 
correspond directly to the remaining tautologies above).

\begin{proposition}\label{prop123}
If $\rho$ contains all the above rules, then
\begin{list}{}{}
\item[(i)]
\; $A \; \turnrho \; B \dimp A$,
\item[(ii)]
\; $A \dimp \falsesent \iffrho -A$,
\item[(iii)]
\; $A \dimp B \iffrho -B \dimp -A$,
\item[(iv)]
\; $A \dimp B, \, B \dimp C \;\, \turnrho \,\; A \dimp C$,
\item[(v)]
\; $A, \, A \dimp B \,\; \turnrho \, \; B$
\item[(vi)]
\; $A \vee B, \, A \dimp C, \, B \dimp C \, \; \turnrho \, \; C$, and
\item[(vii)]
\; if $\Gamma \; \turnrho \; A \dimp B$ then $\Gamma, A \; \turnrho \; B$. 
\end{list}
\end{proposition}

\begin{proof}
$(i)$ Use the Disjunction Introduction Rule.  

$(ii)$ One direction follows from the Disjunction Introduction Rule.  
The other direction uses Proposition~\ref{prop113}$(ii)$  and 
Proposition~\ref{prop96}$(ii)$.

$(iii)$ This follows from Proposition~\ref{prop115}$(i)$, 
Proposition~\ref{prop911},
Proposition~\ref{prop97}$(iii)$, and Proposition~\ref{prop82}$(iii)$.
  
$(iv)$ Use the Disjunction Introduction Rule (twice), 
Proposition~\ref{prop113}$(i)$, and Proposition~\ref{prop96}$(ii)$
(twice).
 
$(v)$ Use Proposition~\ref{prop113}$(i)$ to get $A, -A \; \turnrho \; B$.  
Since $A, B \; \turnrho \; B$, the result follows from Proposition~\ref{prop96}$(ii)$.  

$(vi)$ Use Part$(v)$ to get $A \dimp C, B \dimp C, A \; \turnrho \; C$ 
and  $A \dimp C, B \dimp C, B \; \turnrho \; C$.  
Then the result follows from Proposition~\ref{prop96}$(ii)$. 
 
$(vii)$ This follows from Part$(v)$.
\end{proof}

The last three parts of Proposition~\ref{prop123} show that in some ways
the material conditional~$\dimp$ is stronger than the deductive
conditional~$\imprho$.  
Section~\ref{s14} discusses the corresponding rules obtained by 
replacing $\dimp$ 
with $\imprho$ in the last three parts of Proposition~\ref{prop123}.

The converse of Proposition~\ref{prop123}$(vii)$ does not hold.  
Choose $A$ equal to $B$ where $\ahalt(A)$ is false and $\rho$ is a valid library.  
Then $A \turnrho B$ is true, but $\turnrho A \dimp B$ is false. 
Contrast this with Theorem~\ref{theorem79}.
This illustrates a sense in which~$\dimp$ 
is weaker than~$\imprho$.

Suppose $\ahalt(A)$ and that $\rho$ is valid.
Then $\notrho A$ if and only if $-A$.
Similarly, under these conditions,
$A \imprho B$ if and only if $A \dimp B$.  
The following proposition shows what can be done with a halting 
assumption but without assuming that $\rho$ is valid.

\begin{proposition}\label{prop124}
If  $\rho$ contains all the above rules, then
\begin{list}{}{}
\item[(i)]
\; if $\notrho A$ then $\ahalt(A) \turnrho -A$,
\item[(ii)]
\; $\notrho A \; \turnrho \; \ahalt(A) \imprho -A$,
\item[(iii)]
\; if $\Gamma, A \; \turnrho \; B$ then $\Gamma, \ahalt(A) \; \turnrho \; A \dimp B$,
\item[(iv)]
\; if $A \imprho B$ then $\ahalt(A) \turnrho A \dimp B$, and
\item[(v)]
\; $A \imprho B \; \turnrho \; \ahalt(A) \imprho (A \dimp B)$.
\end{list}
\end{proposition}

\begin{proof}
$(i)$  
By assumption $A \imprho \falsesent$, and $\falsesent \imprho -A$ by Proposition~\ref{prop113}$(ii)$, so $A \turnrho -A$.  
This together with $-A \turnrho -A$ gives the result by Proposition~\ref{prop96}$(ii)$.

$(ii)$ 
This follows from the Universal Rule, Proposition~\ref{prop113}$(ii)$, the
Transitivity rule, and Proposition~\ref{prop910}.

$(iii)$ 
By the Disjunction Introduction Rule, $\Gamma, -A \; \turnrho \; A \dimp B$.  
We have $\Gamma, A \; \turnrho \; A \dimp B$ by assumption
and the Disjunction Introduction Rule.  
The result follows from Proposition~\ref{prop96}$(ii)$.

$(iv)$ This follows from Part$(iii)$.

$(v)$ 
This follows from the Disjunction Introduction, Universal, and Transitivity Rules, 
and Proposition~\ref{prop910}.
\end{proof}

We mentioned above several tautologies of classical logic that do not hold in general 
for the algorithmic material conditional. 
When restricted to algorithmic statements that are true
or directly false, however,  
the algorithmic material conditional can be expected to behave precisely as
the classical material conditional.  
The following corollary illustrates this phenomenon.

\begin{corollary}\label{cor125}
If $\rho$ contain all the above rules, then
\begin{list}{}{}
\item[(i)]
\; $\ahalt(A) \; \turnrho \; A \dimp A$,
\item[(ii)]
\; $\ahalt(A) \; \turnrho \; A \dimp A \vee B$, and
\item[(iii)]
\; $\ahalt(A), B \; \turnrho \; A \dimp A \wedge B$.  
\end{list}
\end{corollary}

\begin{proof}
These follow directly from Proposition~\ref{prop124}$(iv)$ and $(iii)$.
\end{proof}


\section{Stable Base} \label{s13}

The eleven rules developed above are clearly not a complete 
collection of rules for algorithmic logic.  
Indeed, as will be seen in Proposition~\ref{prop137}, one can never have a
complete library of rules for algorithmic logic.  Rather, the rules discussed so far
provides a convenient stable base on which to build more elaborate
stable libraries.

\begin{definition} \label{def131}
A \emph{base} is a set $\cB$ of rules.  
We require that a base be finite, or at least 
arises as the set of terms of a library.
A \emph{$\cB$-library} is a library containing all the rules
of the base $\cB$.  
A $\cB$-library~$\rho$ is said to be \emph{valid outside~$\cB$}
if every rule in~$\rho$ which is not in~$\cB$ is $\rho$-valid.
A base $\cB$ is \emph{stable} if every $\cB$-library
that is valid outside of~$\cB$ is itself valid.

Let $\cB_0$ be the set containing
Rules~\ref{rule1} to \ref{rule11} above.
\end{definition}

\begin{theorem} \label{theorem132}
The set $\cB_0$ is a stable base.
\end{theorem}

\begin{proof}
Let $\rho$ be a $\cB_0$-library that is valid outside $\cB_0$.
We need to show that $\rho$ is valid.

Rules~\ref{rule1}, \ref{rule2}, \ref{rule4}, \ref{rule6},
and \ref{rule9} are $\rho$-valid by
Propositions~\ref{prop51}, \ref{prop61}, \ref{prop72}, \ref{prop92},
\ref{prop112}, respectively.  
Rules~\ref{rule10} and~\ref{rule11} are $\rho$-valid
by Proposition~\ref{prop114}.
Rule~\ref{rule3} is $\rho$-valid by Proposition~\ref{prop64} since
$\cB$ contains the Universal Rule.
Rule~\ref{rule5} is $\rho$-valid by Proposition~\ref{prop76}
since $\cB$ contains the Conjunction Rule.
Rule~\ref{rule8} is $\rho$-valid by Proposition~\ref{prop99} since $\cB$
contains the Disjunction Elimination, Universal, and Conjunction
Rules.  
Finally, Theorem~\ref{theorem95} takes care of Rule~\ref{rule7}
and shows that $\rho$ is valid.
\end{proof}

\begin{definition} \label{def133}
Let $\cB$ be a base.  A rule is \emph{$\cB$-safe} if it is $\rho$-valid for all
$\cB$-libraries~$\rho$.  A \emph{stable extension} $\cB'$ of $\cB$ is a base
containing $\cB$ such that every rule in $\cB'$ that is not in $\cB$ is $\cB$-safe.
\end{definition}

\begin{proposition}\label{prop134}
A stable extension of a stable base is a stable base.
\end{proposition}

\begin{proof}
Let $\cB$ be a stable base and $\cB'$ a stable extension of $\cB$.    
Suppose $\rho$ is a $\cB'$-library
valid outside of $\cB'$. We must show that $\rho$ is valid. 

First we show that $\rho$ is actually valid outside~$\cB$.
To that end, let $\rho_k$ be outside $\cB$. If $\rho_k$
happens to be in~$\cB'$ then it is $\cB$-safe by the definition of
stable extension. In particular,  $\rho_k$ is $\rho$-valid.
If $\rho_k$ is outside $\cB'$ then it
is $\rho$-valid simply because $\rho$ is valid outside of~$\cB'$.
Thus $\rho$ is valid outside $\cB$. 

Since $\cB$ is stable, and since $\rho$ is valid outside~$\cB$,
the library $\rho$ is valid.
\end{proof}

\begin{proposition}\label{prop135}
If the rules of a library forms a stable base, then the library is valid.
Thus, if the rules of a library form a stable extension of $\cB_0$,
then the library is valid.
\end{proposition}

\begin{definition}\label{def136}
Let $\rho_1$ and $\rho_2$ be libraries. 
Then $\rho_2$ is \emph{stronger than} $\rho_1$ if
$A \imprhoone B$ implies $A \imprhotwo B$ for all
algorithmic statements $A$ and $B$.  A library
$\rho_2$ is \emph{strictly stronger than} $\rho_1$ if
$(i)$ $\rho_2$ is stronger than $\rho_1$, and $(ii)$
there exists
$A$ and $B$ such that $A \imprhotwo B$ is true but
$A \imprhoone B$ is false.
\end{definition}

Algorithmic logic is complete in the sense that, for any library~$\rho$
containing the Universal Rule,  if $A$ is true
then $\proverho(A)$ is true. But there is also a sense in which the
logic is inherently incomplete.

\begin{proposition}\label{prop137}
For every valid library $\rho_1$ there is a strictly stronger valid library $\rho_2$. 
If the set of rules in $\rho_1$ form a stable base~$\cB$, 
then $\rho_2$ can be taken to be a library whose rules
form a stable $\cB$-extension.
\end{proposition}

\begin{proof}
Let $\mpone$ be the algorithm implementing 
$\begin{array}{c}
A \imprhoone B\\
A\\
\cline{1-1}
B
\end{array}
$.
This is a $\rho$-valid rule for any library~$\rho$
since $\rho_1$ is valid.\footnote{
This rule differs essentially from $P_3$ discussed
in Section~\ref{s14} in that, 
given input $[H, \rho, m]$, the rule $\mpone$ does not use the
input~$\rho$, but rather uses
the fixed library~$\rho_1$.  Rule $P_3$, on the other hand, does
use the input $\rho$.}
There is no algorithm that decides whether a statement is false.
Thus there is a false algorithmic statement $C$ such
that $C \imprhoone \falsesent$ is false.  
Let $\denyc$ be the algorithm implementing 
$\begin{array}{c}
C\\
\cline{1-1}
\falsesent
\end{array} 
$.
The rule $\denyc$ is $\rho$-valid for any library $\rho$ since
$C$ is false.

Let $\rho_2$ be the library containing $\mpone$, $\denyc$,
and the Universal Rule.
Observe that $(i)$ $\rho_2$ is valid, $(ii)$ if $A \imprhoone B$
then $A \imprhotwo B$, and $(iii)$ $C \imprhoone \falsesent$ is
false but $C \imprhotwo \falsesent$ is true.  To see $(ii)$, let $S$
be the $\rho_2$-deductive closure of $A$.  By the Universal Rule,
$A \imprhoone B$ is in $S$.  By $\mpone$, $B$ is in $S$.  So
$\rho_2$ is valid and strictly stronger than $\rho_1$.

Now suppose that the rules of $\rho_1$ form a stable base $\cB$. 
Let $\rho_2$ contain all the rules of $\rho_1$, $\mpone$, $\denyc$,
and the Universal Rule. The new rules are $\cB$-safe, so
the rules of $\rho_2$ form a stable valid extension of $\cB$.  Finally,
by an argument similar to the one above, $\rho_2$ is strictly
stronger than $\rho_1$.

\end{proof}


\section{Paradoxical Rules}\label{s14}

A \emph{paradoxical rule} is an algorithmic counterpart of a traditional
rule of logic that cannot be in any stable base.\footnote{The 
term \emph{paradoxical} is used since many of the
arguments related to such rules are akin to those occurring in
the Russell and Curry paradoxes.}
In this section we will show that the following are paradoxical rules:

$$
\begin{array}{rrcrcccccc}
P_1\!:\!\!\!\!\!&&
\quad\quad P_2\!:\!\!\!\!\!&&
\quad\quad P_3\!:\!\!\!\!&&
\quad P_4\!:&A \imprho C&
\quad P_5\!:\!\!\!\!
\\
&A&&A&&A\imprho B&&B \imprho C&&\\
&\notrho A&&\notrho A&& A&&A \vee B&&\notrho\; \notrho A\\
\cline{2-2}\cline{4-4}\cline{6-6}\cline{8-8}\cline{10-10}
&\falsesent && B && B && C && A \\
\end{array}
$$

$$
\begin{array}{cccccccccc}
P_6\!:\!\!\!&A \vee B&
\;\; P_7\!:\!\!\!&\notrho A \vee B&
\;\;\; P_8\!:\!\!\!\!\!\!\!\!&&
\;\;\; P_9\!:\!\!\!\!\!\!\!\!&&
P_{10}\!:\!\!\!\!&
\\
&\notrho A&&A&&\emptyset && \emptyset&&A\imprho B\\
\cline{2-2}\cline{4-4}\cline{6-6}\cline{8-8}\cline{10-10}
&B&&B&&
A \vee \notrho A&\quad& \notrho A \vee \notrho\; \notrho A&&
\notrho A \vee B
\\
\end{array}
$$

$$
\begin{array}{cccccccc}
P_{11}\!:\!\!\!\!&&
\; P_{12}\!:\!\!\!\!\!\!\!&&
\;P_{13}\!:\!\!\!\!\!&&
P_{14}\!:\!\!\!\!&\\
&\notrho (A \wedge B)&& A && \proverho (A) && \proverho (\proverho (A))\\
\cline{2-2}\cline{4-4}\cline{6-6}\cline{8-8}
&\notrho A \vee \notrho B  &&  \notrho\; \notrho A &&  A &&  \proverho(A)\\
\end{array}
$$
Given the expected input $[H, \rho, m]$, all of the rules above use
$\rho$, but only Rules $P_2$, $P_8$, and $P_9$ use the
resource integer $m$ in their implementation.  The symbol $\emptyset$
in Rules $P_8$ and $P_9$ indicates that no premises in $H$ are required.  
Clearly, some of
the paradoxical rules above
are interrelated.

Rules $P_1$, $P_3$, $P_4$, $P_6$, $P_7$, $P_{13}$,
and $P_{14}$ have the remarkable property of being $\rho$-valid
for any valid library~$\rho$ but, due to their instability,
not being in any sufficiently rich valid~$\rho$.
Rule $P_2$ has a similar status, at least for any $\cB_0$-library~$\rho$.

\begin{remark}
As one might expect, many of these correspond to rules that have aroused
suspicion in the past and have been excluded from
weaker logics such as intuitionistic or minimal logic.
The long list of paradoxical rules 
to be avoided in algorithmic logic might make
algorithmic logic seem weak. However, 
in algorithmic logic one always has the option of going to a stronger library $\rho$,
often compensating for not having the above rules.
\end{remark}

\begin{lemma} \label{lemma141} 
Every stable base $\cB$ has a stable extension $\cB'$
with the following property:  For every $\cB'$-library $\rho$
there is an algorithmic statement~$Q_\rho$ such that
$Q_\rho \iffrho \notrho Q_\rho$. 
\end{lemma}

\begin{proof}
The proof requires an algorithm $\curry$ that expects as input
a list $[\alpha, \rho]$ where $\alpha$ is an algorithm.  
If the algorithmic statement $\notrho \bigl[\alpha, [\alpha, \rho], 1\bigr]$ is true, 
then $\curry$ outputs $1$.  
Otherwise $\curry$ does not halt.\footnote{This
algorithm is called $\curry$ due to the resemblance
of the proof of Theorem~\ref{theorem142} to a common version
of the Curry Paradox.}
Observe that if $\alpha$ is an algorithm, 
then $\bigl[ \curry, [\alpha, \rho], 1 \bigr]$ if and only if
$\notrho \bigl[ \alpha, [\alpha, \rho], 1 \bigr]$.

Let $\beta$ be the rule that simultaneously implements the two rule diagrams:
$$
\begin{array}{cccl}
\bigl[\curry, [\alpha, \rho], 1\bigr]&&
\notrho \bigl[\alpha, [\alpha,\rho], 1\bigr]\\
\cline{1-1}\cline{3-3}
\notrho \bigl[\alpha, [\alpha, \rho], 1\bigr] && 
\bigl[\curry, [\alpha, \rho], 1\bigr]&.
\end{array}
$$
More specifically, assuming an input of the expected form $[H, \rho, m]$,
the rule $\beta$ looks for all statements of the form of the first
line of either of the above diagrams, where $\alpha$ is required
to be an algorithm. For each such statement it finds, it appends
the appropriate statement to~$H$.

Clearly $\beta$ is $\cB$-safe where $\cB$ is the given stable base.
Let $\cB'$ be the stable extension of~$\cB$ obtained by simply adding
the rule~$\beta$ to $\cB$.  
Given a $\cB'$-library~$\rho$, let
$Q_\rho$ be $\bigl[\curry, [\curry, \rho], 1\bigr]$. 
So $Q_\rho \iffrho \notrho Q_\rho$ since
$\rho$ contains $\beta$.
\end{proof}

While the rule $\beta$ used in the above proof is not a rule of 
elementary logic, and may thus seem \emph{ad hoc}, 
it is a consequence of general, more natural rules concerning
the basic properties of algorithms.
This is discussed in~\cite{thirdpaper06}.

\begin{theorem} \label{theorem142}
There is no stable base $\cB$ such that 
the law $A, \notrho A \; \turnrho \; \falsesent$
holds for all valid $\cB$-libraries $\rho$.
\end{theorem}

\begin{proof}
Suppose that there is such a $\cB$, and
let $\cB'$ be as in Lemma~\ref{lemma141}.  
Let $\rho$ be the library consisting of the 
rules of $\cB'$. The validity of $\rho$
follows from Proposition~\ref{prop135}.
By Lemma~\ref{lemma141}, 
there is a statement $Q_\rho$ such that
$Q_\rho \iffrho \notrho Q_\rho$. By validity, 
$Q_\rho$ holds if and only if $\notrho Q_\rho$
holds. 

Let $S$ be the $\rho$-deductive closure of $Q_\rho$.  
Since, $Q_\rho \, \turnrho\!\! \notrho Q_\rho$, the set $S$ 
contains~$\notrho Q_\rho$.  
By assumption $Q_\rho, \notrho Q_\rho \; \turnrho \; \falsesent$
holds, so $S$ contains $\falsesent$.  
Thus $Q_\rho \imprho \falsesent$ holds; that is, $\notrho Q_\rho$ is true.  
As mentioned above, this implies that $Q_\rho$ is true.
Since $Q_\rho, \notrho Q_\rho \; \turnrho \; \falsesent$
and since $\rho$ is valid, $\falsesent$ is true. 
\end{proof}

\begin{corollary} \label{cor143}
No stable base contains Rule $P_1$ or Rule $P_2$
(defined at the beginning of this section).
\end{corollary}

\begin{corollary} \label{cor144}
If $\cB$ is a stable base, then the assertion that
$$
\hbox{$\Gamma \; \turnrho \; A \imprho B$ implies
$\Gamma, A \; \turnrho\; B$}
$$ 
fails for some valid $\cB$-library $\rho$.
\end{corollary}

\begin{proof}
Let $\rho$ be a valid $\cB$-library for which the assertion holds.
By Corollary~\ref{cor415}(ii), $\notrho A \; \turnrho \; \notrho A$.
In other words, $\notrho A \; \turnrho \; A\imprho \falsesent$.
So $\notrho A, A  \; \turnrho \; \falsesent$ by the assertion.
By Theorem~\ref{theorem142} this cannot hold for all such
$\rho$.
\end{proof}

\begin{corollary} \label{cor145}
There is no stable base $\cB$ such that 
$A \imprho B, A \; \turnrho \; B$
holds for all $\cB$-libraries $\rho$.  
In particular, no stable base contains Rule~$P_3$.
\end{corollary}

\begin{proof}
Suppose otherwise. If $\rho$ is a valid $\cB$-library,
then $A, A \imprho \falsesent \; \turnrho \; \falsesent$
for all~$A$.  In other
words, $A, \notrho A \; \turnrho \; \falsesent$ holds,
contradicting Theorem~\ref{theorem142}.
\end{proof}

\begin{corollary}\label{cor146}
There is no stable base $\cB$ such that the law
$$A \imprho C, \; B \imprho C, \; A \vee B \;\; \turnrho \;\; C$$
holds for all $\cB$-libraries $\rho$. 
In particular, no stable base contains Rule $P_4$. 
\end{corollary}

\begin{proof}
Suppose that there is such a $\cB$. 
The Disjunction Introduction Rule is $\cB$-safe, so the base $\cB'$ that
results from adding this rule to $\cB$ is also stable.
Let $\rho$ be a valid $\cB'$-library.  
Let $S$ be the deductive closure of $A$ and $\notrho A$.
By the Disjunction Introduction Rule, $A \vee A$ is in $S$.  Since
$A \imprho \falsesent$ is in $S$, so is $\falsesent$.
Thus $A, \notrho A \; \turnrho \; \falsesent$ for all $A$ and 
all such~$\rho$, contradicting
Theorem~\ref{theorem142} for the stable base~$\cB'$.
\end{proof}

\begin{corollary}\label{cor147}
There is no stable base $\cB$ such that the 
law $\notrho\; \notrho A \; \turnrho \; A$
holds for all $\cB$-libraries $\rho$.  
In particular, no stable base contains Rule~$P_5$.
\end{corollary}

\begin{proof}
Suppose that there is such a stable base $\cB$. 
The Universal and Transitivity Rules are $\cB$-safe, 
so the base $\cB'$ that results from adding these rules to $\cB$ is also stable.  
The Meta-Universal Rule is $\cB'$-safe since $\cB'$ contains the Universal Rule, 
so the base $\cB''$ that results from adding the Meta-Universal Rule to $\cB'$ 
is also stable.

Let $\rho$ be a valid $\cB''$-library and $A$ a statement.  
Let $S$ be the deductive
closure of $A$ and $\notrho A$.
By the Meta-Universal rule $\notrho \falsesent \imprho A$ is in~$S$.  
Since $\notrho A$ is $A \imprho \falsesent$, which is in~$S$, 
$\notrho \falsesent \imprho \falsesent$ is in~$S$ by the Transitivity Rule.
In other words, $\notrho\; \notrho \falsesent$ is in~$S$.
Thus $\falsesent$ is in~$S$ by hypothesis.
So $A, \notrho A \; \turnrho \; \falsesent$,
contradicting Theorem~\ref{theorem142} for the base $\cB''$.
\end{proof}

\begin{corollary} \label{cor148}
There is no stable base $\cB$ where the law
$A \vee B, \notrho A \; \turnrho \; B$ holds for all
$\cB$-libraries $\rho$.  
Likewise, there is no stable base $\cB$ where the law
$\notrho A \vee B, A \; \turnrho \; B$ holds for all
$\cB$-libraries $\rho$.  
In particular, no stable base contains either Rule
$P_6$ or Rule $P_7$.
\end{corollary}

\begin{proof}
Suppose that there is a $\cB$ where the first of these
laws holds.  The Disjunction Introduction
Rule is $\cB$-safe, so the extension $\cB'$ obtained by adding 
this Rule to~$\cB$ is stable.
Let $\rho$ be any valid $\cB'$-library.

Let $S$ be the deductive closure of $A$ and $\notrho A$.
By the Disjunction Introduction Rule, $A \vee \falsesent$ is in $S$.  
By hypothesis, $\falsesent$ is in~$S$.
We have established that $\notrho A, A \; \turnrho \falsesent$ holds
for every statement~$A$ and valid $\cB'$-library~$\rho$, contradicting
Theorem~\ref{theorem142}. 

The second part of the theorem follows by a similar argument.
\end{proof}

\begin{lemma}\label{lemma149}
Suppose $A$ is an algorithmic statement where $A \iffrho \notrho A$ 
with $\rho$ a valid library.  
Then $A$, $\notrho A$, and $\notrho \; \notrho A$ are all false.
\end{lemma}

\begin{proof}
Suppose $A$ is true.  
By hypothesis, $A \imprho \; \notrho A$.  
So, by Proposition~\ref{prop418}$(iii)$ and the validity of $\rho$, 
the statement $\notrho A$ holds.  
Thus $A$ and $A \imprho \falsesent$ hold. 
Again, by Proposition~\ref{prop418}$(iii)$, $\falsesent$ is true.

Suppose $\notrho A$.  By hypothesis, $ \notrho A \imprho A$.  So
$A$ is true, contradicting the above.

Suppose $\notrho \; \notrho A$; in other words, $\notrho A \imprho \falsesent$.
By hypothesis, $A \imprho \notrho A$. 
By Proposition~\ref{prop418}$(ii)$, $A \imprho \falsesent$.
In other words, $\notrho A$ which contradicts the above.
\end{proof}

\begin{corollary}\label{cor1410} 
Let $\cB$ be a stable base.  There is a valid $\cB$-library
$\rho$ and an algorithmic statement $Q_\rho$ such that
$Q_\rho \vee \notrho Q_\rho$ and
$\notrho Q_\rho \vee \notrho \; \notrho Q_\rho$ are both false.  In
particular, Rules $P_8$ and $P_9$ are not $\rho$-valid. 

Therefore, there is no stable base containing Rules $P_8$ or $P_9$.
\end{corollary}

\begin{proof}
Let $\cB'$ be as in Lemma~\ref{lemma141}.
Let $\rho$ be a library consisting of the rules in~$\cB'$, 
and let  $Q_\rho$ be as in Lemma~\ref{lemma141}.
The result follows from Lemma~\ref{lemma149}.
\end{proof}

\begin{proposition} \label{prop1411}
There is no stable base $\cB$ where
$A \imprho B \, \turnrho  \notrho A \vee B$ holds for all
$\cB$-libraries $\rho$.  
In particular, no stable base contains Rule~$P_{10}$.
\end{proposition}

\begin{proof}
Suppose that there is such a~$\cB$.  
Let $\rho$ be a library consisting of the rules in
$\cB'$ as defined in Lemma~\ref{lemma141}.  
Let $Q_\rho$ be as in Lemma~\ref{lemma141}.
The library $\rho$ is valid since $\cB'$ is a stable base.  

By hypothesis
$Q_\rho \imprho \notrho Q_\rho \; \turnrho \; \notrho Q_\rho \vee \notrho Q_\rho$, 
so by Proposition~\ref{prop416} and the validity of $\rho$, 
the statement $\notrho Q_\rho \vee \notrho Q_\rho$ is true.  
So $\notrho Q_\rho$ is true contradicting Lemma~\ref{lemma149}.
\end{proof}

\begin{theorem} \label{theorem1412}
There is no stable base $\cB$ where the law
$$\notrho (A \wedge B) \; \turnrho \; \notrho A \vee \notrho B$$ 
holds for all $\cB$-libraries~$\rho$.  
In particular, no stable base contains Rule~$P_{11}$.
\end{theorem}

\begin{proof}
Suppose that there is such a stable base $\cB$.
The rule represented by the diagram
$ \begin{array}{c}
A \wedge -A \\
\cline{1-1}
\falsesent\\
\end{array}
$
is $\cB$-safe.  
Let $\cB'$ be the stable extension obtained by adding
this rule to $\cB$.  
Let $\rho$ be a library consisting of the rules in~$\cB'$.  
The library $\rho$ is valid since $\cB'$ is a stable base.
The statement $\notrho (A \wedge -A)$ holds for any $A$ because of
the new rule added to the library.  
So, by hypothesis and the validity of $\rho$,
the statement $\notrho A \vee \notrho - A$ is true for all $A$.

Let $\beta$ be an algorithm that expects an algorithm $\alpha$ as input.  
The algorithm $\beta$ finds the smallest $m$ such
that $\notrho [\alpha, \alpha, 1]$ 
or $\notrho - [\alpha, \alpha, 1]$ is $m$-true.
There will be such an $m$ since
 $\notrho A \vee \notrho - A$ holds for all $A$.
If $\notrho [\alpha, \alpha, 1]$ is $m$-true for
this value of $m$, then $\beta$ outputs $1$.  
If  $\notrho [\alpha, \alpha, 1]$ is not $m$-true, but $\notrho - [\alpha, \alpha, 1]$ 
is $m$-true for this value of $m$, then $\beta$ outputs~$0$.  
The notion of $m$-true here is as in the definition of the Universal Rule.
Observe that if $\alpha$ is an algorithm then $\beta$ halts 
for input~$\alpha$.

Let $B$ be the statement $[\beta, \beta, 1]$.  
If $B$ is true, then $\notrho B$ is true.
If $-B$ is true, then $\notrho -B$.  
Since $\rho$ is valid, it is not possible for a statement $A$ and its 
negation $\notrho A$ to both be true.  
So neither $B$ nor $-B$ is true.
In other words, $\beta$ does not halt for input~$\beta$, a contradiction.
\end{proof}

\begin{proposition}\label{prop1413}
Let $\cB$ be a stable base.  
Then there is a valid $\cB$-library $\rho$ 
such that $\notrho \; \notrho \truesent$ is false.  
Furthermore, the law $A \; \turnrho \; \notrho \; \notrho A$ does 
not hold for all $\cB$-libraries $\rho$.  
In particular, no stable base contains Rule~$P_{12}$.
\end{proposition}

\begin{proof}
As in the proof of Corollary~\ref{cor147}, there is a stable
base $\cB''$ containing $\cB$ together with the Universal, the
Meta-Universal, and the Transitivity Rules.  
Let $\rho$ be any valid $\cB''$-library.
Suppose, $\notrho \; \notrho \truesent$ holds.
Thus $\notrho \truesent \turnrho \falsesent$.
Let $S$ be the deductive closure of $A$ and $\notrho A$.
By the Meta-Universal Rule, $\truesent \imprho A$ is in~$S$. 
Note that $A \imprho \falsesent$ is in $S$, so, by the Transitivity Rule,
$\truesent \imprho \falsesent$ is in~$S$.  
In other words, $\notrho \truesent$ is in~$S$.  
Since $\notrho \truesent \turnrho \falsesent$ is true, 
$\falsesent$ must be in~$S$.

We have established that if 
$\notrho \; \notrho \truesent$ holds
then $\notrho A, A \; \turnrho \falsesent$
holds for all~$A$.
Therefore, by Theorem~\ref{theorem142}, there must be a valid
$\cB''$-library~$\rho$ such that  $\notrho \; \notrho \truesent$ is false.
The law $A \; \turnrho \; \notrho \; \notrho A$ does not hold for such~$\rho$.  
To see this, consider the case where $A$ is $\truesent$.
\end{proof}

\begin{theorem} \label{theorem1414} 
There is no stable base $\cB$ where the law
$$\proverho (\proverho (A)) \; \turnrho \; \proverho(A)$$
holds for all $\cB$-libraries $\rho$.  
In particular, there is no stable base $\cB$ where
$\proverho(A) \; \turnrho \; A$ holds for all $\cB$-libraries $\rho$.
So no stable base contains
Rule~$P_{13}$ or Rule~$P_{14}$.
\end{theorem}

\begin{proof}
Suppose otherwise that there is such a stable base $\cB$.
As in the proof of Corollary~\ref{cor147}, there is a stable
base $\cB''$ containing $\cB$ together with the Universal, the
Meta-Universal, and the Transitivity Rules. 

Consider an algorithm~$\beta$ that expects
as input~$[\alpha, \rho]$ 
where $\alpha$ is an algorithm. The algorithm
$\beta$ checks the truth of
$\bigl[\alpha, [\alpha, \rho], 1\bigr] \imprho \proverho(\falsesent)$.
If the statement is true, $\beta$ outputs $1$. Otherwise,
$\beta$ does not halt.
Let $R_\rho$ be the statement $\bigl[\beta, [\beta, \rho], 1\bigr]$.  Observe that
$R_\rho$ is true if and only if $R_\rho \imprho \proverho(\falsesent)$
is true.

The rule implementing
$ \begin{array}{cl}
R_\rho\\
\cline{1-1}
R_\rho \imprho \proverho(\falsesent)
\end{array}
$
is $\rho$-valid for all libraries~$\rho$. 
Let  $\rho$ be the library consisting of this rule together
with all the rules of~$\cB''$.
The $\cB''$-library $\rho$ is valid since $\cB''$ is stable.

Let $S$ be the deductive closure of $R_\rho$.  So
$R_\rho \imprho \proverho(\falsesent)$ is in $S$.
By Proposition~\ref{prop65}$(ii)$,
$\proverho (\proverho(\falsesent))$ is in $S$.
Finally, by supposition, $\proverho (\falsesent)$ is also in $S$.
We have shown that $R_\rho \imprho \proverho (\falsesent)$
is true.  Therefore, $R_\rho$ is true.  Since $\rho$ is valid,
Proposition~\ref{prop418}$(iii)$ implies that $\proverho( \falsesent)$
is true.  So by Proposition~\ref{prop418}$(iv)$ and the validity of $\rho$,
$\falsesent$ is true.
\end{proof}


\section{Conclusion}\label{s15}

In~\cite{thirdpaper06} we introduce additional rules to algorithmic
logic which do not concern logical connectives as do the rules in the
current paper. Instead, these new rules
relate to the basic structure of algorithms themselves.
These structural rules will lead to a strong internal abstraction
principle making algorithmic logic more flexible and powerful. 

In particular, for bases~$\cB$ containing these structural rules, 
Lemma~\ref{lemma141} can be strengthened to apply 
to all $\cB$-libraries~$\rho$. 
Consequently, the main results of Section~\ref{s14} can be
significantly strengthened. 
More precisely, let $\cB_1$ be the stable base consisting of $\cB_0$
together with the structural rules of the promised future paper.
Many of the results of Section~\ref{s14} refer to laws which do not hold
for all $\cB$-libraries. In other words
there exists \emph{some} $\cB$-library where the law fails.
For the base $\cB_1$, however,  these results can be strengthened to assert that
the given law fails for \emph{all} valid $\cB_1$-libraries.

In Section~\ref{s14} above we mention that several of the paradoxical rules are
$\rho$-valid as long as $\rho$ is valid. 
The other rules, with one exception, cannot be expected to be $\rho$-valid.
More specifically, if $\rho$ is a valid $\cB_1$-library, then
all the other rules, with the exception of Rule~$P_5$, are not $\rho$-valid.
This can be seen with arguments similar to those of Section~\ref{s14}.
Rule~$P_{5}$ is $\rho$-valid for such~$\rho$, however, because of the striking fact
that $\notrho\; \notrho A$ is false for \emph{all}~$A$.
This fact can be shown with an argument 
similar to that of Proposition~\ref{prop1413}.\footnote{We would like to thank
our colleagues for many useful discussions, and the referees for
several good suggestions including 
the suggestion to use the term \emph{strong negation} in honor of David
Nelson. 
One referee asked an interesting question concerning
the status of~$(A\imprho \notrho B)\; \turnrho \; (B  \imprho \notrho A)$,
a contrapositive law whose analogue holds in intuitionistic and even minimal logic.
For sufficiently rich $\rho$, the statement $\notrho\; \notrho \truesent$
is false (Proposition~\ref{prop1413}).
For such valid $\rho$ the law fails: consider the case where
$A$ is $\notrho \truesent$ and $B$ is $\truesent$.}



\end{document}